\documentclass{amsart}
\usepackage{amscd,amsmath,amssymb,euscript}

\newtheorem{theorem}{Theorem}[section]
\newtheorem{proposition}[theorem]{Proposition}
\newtheorem{lemma}[theorem]{Lemma}
\newtheorem{definition}[theorem]{Definition}
\newtheorem{remark}[theorem]{Remark}
\newtheorem{example}[theorem]{Example}
\newcommand{\Ext}{\mathop\mathrm{Ext}\nolimits}
\newcommand{\Hom}{\mathop\mathrm{Hom}\nolimits}

\newcommand{\coker}{\mathop\mathrm{coker}\nolimits}
\newcommand{\im}{\mathop\mathrm{im}\nolimits}
\newcommand{\rk}{\mathop\mathrm{rank}\nolimits}
\newcommand{\Spec}{\mathop\mathrm{Spec}\nolimits}
\newcommand{\Quot}{\mathop\mathrm{Quot}\nolimits}

\newcommand{\sch}{\mathrm{Sch}}
\newcommand{\sets}{\mathrm{Sets}}

\begin{document}

\title{Moduli of stable objects
in a triangulated category}
\author{Michi-aki Inaba}
\address{{\rm Michi-aki Inaba} \\
Department of Mathematics, Kyoto University, \\
Kyoto, 606-8502, Japan}
\email{inaba@math.kyoto-u.ac.jp}
\footnote{Key Words and Prases: moduli, triangulated category}

\maketitle

\begin{abstract}
We introduce the concept of strict ample sequence in 
a fibered triangulated category and define the stability
of the objects in a triangulated category.
Then we construct the moduli space of (semi) stable objects
by GIT construction.
\end{abstract}

\section{Introduction}

Let $X\rightarrow S$ be a projective and flat morphism
of noetherian schemes.
We consider the functor
$\mathrm{Splcpx}_{X/S}:(\sch/S)\rightarrow(\sets)$
defined by
\[
 \mathrm{Splcpx}_{X/S}(T)=\left\{
 E\in D^b(\mathrm{Coh}(X\times_ST)) \left|
 \begin{array}{l}
  \text{for any geometric point $t$ of $T$,} \\
  \text{$E(t):=E\otimes^L k(t)$ is a bounded complex and} \\
  \Ext^i(E(t),E(t))\cong \left\{
   \begin{array}{l}
    k(t) \quad \text{if $i=0$} \\
    0 \quad \text{if $i=-1$}
   \end{array} \right.
 \end{array}
 \right\}\right/\sim,
\]
where $E\sim E'$ if there is a line bundle $L$ on $T$ such that
$E\cong E'\otimes L$ in $D^b(\mathrm{Coh}(X\times_ST))$.
We denote the \'{e}tale sheafification of $\mathrm{Splcpx}_{X/S}$
by $\mathrm{Splcpx}_{X/S}^{\text{\'{e}t}}$.
Then the result of \cite{inaba} is that
$\mathrm{Splcpx}_{X/S}^{\text{\'{e}t}}$
is an algebraic space over $S$.
M. Lieblich extends this result in \cite{Lieblich}
to the case when $X\rightarrow S$ is a proper flat morphism
of algebraic spaces.
So the problem on the construction of the moduli space of objects in 
a derived category is solved in some sense.
However, the moduli space $\mathrm{Splcpx}_{X/S}^{\text{\'et}}$
is not separated and it is not a good space in geometric sense.
So we want to construct a projective moduli space
(or quasi-projective moduli space with a good compactification)
as a Zariski open set of $\mathrm{Splcpx}_{X/S}^{\text{\'et}}$
such as the moduli space of stable sheaves.

This problem is also motivated by Fourier-Mukai transform.
Let $X,Y$ be projective varieties over an algebraically closed field $k$
and ${\mathcal P}$ be an object of $D^b(\mathrm{Coh}(X\times Y))$.
The functor
\begin{gather*}
 \Phi: D^b(\mathrm{Coh}(X)) \longrightarrow D^b(\mathrm{Coh}(Y)) \\
  \hspace{30pt} E \mapsto \mathrm{R}(p_Y)_*(p_X^*(E)\otimes^L{\mathcal P})
\end{gather*}
is called a Fourier-Mukai transform if it is an equivalence of categories.
Here $p_X:X\times Y\rightarrow X$ and $p_Y:X\times Y\rightarrow Y$
are the projections.
Fourier-Mukai transform induces the isomorphisms on moduli spaces and
for example the image $\Phi(M_X^P)$ of a moduli space of stable sheaves
$M_X^P$ on $X$ by $\Phi$ sometimes becomes a moduli space of stable sheaves
on $Y$.
The problem on the preservation of stability under Fourier-Mukai transform
is investigated by many people and this problem is clearly pointed out by
 K. Yoshioka in \cite{yoshioka2}.
However, the image $\Phi(M_X^P)$ of the moduli space of stable sheaves
by the Fourier-Mukai transform may not be contained in the category
of coherent sheaves on $Y$ in general and so we must consider
certain moduli space of stable objects in the derived category
$D^b(\mathrm{Coh}(Y))$.

In this paper we introduce the concept ``strict ample sequence"
in a triangulated category.
``Strict ample sequence" satisfies the condition of ample sequence
defined by A. Bondal and D. Orlov in \cite{B-O}, but it also satisfies
many other conditions because we expect that a ``polarization"
is determined by strict ample sequence.
Indeed we can define stable objects determined by a strict ample sequence
and construct the moduli space of stable objects (resp.\ $S$-equivalences
classes of semistable objects) as a quasi-projective scheme
(resp.\ projective scheme).
This is the main result of this paper (Theorem \ref{thm:moduli-exists}
and Theorem \ref{projective}).
If $\Phi:D^b(\mathrm{Coh}(X))\rightarrow D^b(\mathrm{Coh}(Y))$
is a Fourier-Mukai transform and $M_X^P$ is a moduli space of stable
sheaves on $X$, then the image $\Phi(M_X^P)$ of $M_X^P$ by $\Phi$ becomes
a moduli space of stable objects in $D^b(\mathrm{Coh}(Y))$
whose stability is determined by some strict ample sequence on
$D^b(\mathrm{Coh}(Y))$.
So Fourier-Mukai transform always preserves certain stability in our sense
(Example \ref{ex:F-M}).

T. Bridgeland defined in \cite{Bridgeland1} the concept of stability condition
on a triangulated category.
So we are interested in the relation between the stability condition of
Bridgeland and the deifinition of stability determined by a strict
ample sequence.
However, it seems rather impossible to expect the construction of a strict
ample sequence from the stability condition defined by Bridgeland
without any other condition.
How to treat the relation between strict ample sequence and
stability condition of Brigeland is a problem still unsolved.

\section{Definition of fibered triangulated category}

Let $S$ be a noetherian scheme.
We denote the category of noetherian schemes over $S$ by $(\sch/S)$
and the derived category of bounded complexes of coherent sheaves on $U$
by $D^b_c(U)$ for $U\in(\sch/S)$.
We denote the derived category of lower bounded complexes of
coherent sheaves on $U$ by $D^+_c(U)$ for $U\in(\sch/S)$.
For a noetherian scheme $X$ over $S$,
we denote the full subcategory of
$D^b_c(X)$ consisting of the objects
of finite Tor-dimension over $S$
by $D^b(\mathrm{Coh}(X/S))$.
Then $D^b(\mathrm{Coh}(X/S))$ becomes
a triangulated category.
For a triangulated category ${\mathcal T}$ and for objects $E,F\in{\mathcal T}$,
we write $\Ext^i(E,F):=\Hom_{\mathcal T}(E,F[i])$.


\begin{definition}\rm
$p:{\mathcal D}\rightarrow (\sch/S)$
is called a fibered triangulated category if
\begin{enumerate}
\item ${\mathcal D}$ is a category, $p$ is a covariant functor,
\item for any $U\in(\sch/S)$, the full subcategory
 ${\mathcal D}_U:=p^{-1}(U)$ of ${\mathcal D}$ is a triangulated category,
\item for any object $E\in{\mathcal D}_U$ and for any morphism
$f:V\rightarrow U=p(E)$ in $(\sch/S)$,
there exist an object $F\in{\mathcal D}_V$
and a morphism $u:F\rightarrow E$
satisfying the condition:
For any object $G\in{\mathcal D}_V$ and a morphism $v:G\rightarrow E$
with $p(v)=f$, there exists a unique morphism $w:G\rightarrow F$
satisyfing $p(w)=\mathrm{id}_V$ and $u\circ w=v$,
(we denote $F$ by $f^*(E)$ or $E_V$ and we call such morphism $u$ a Cartesian
morphism),
\item any composition of Cartesian morphisms is Cartesian,
\item for any morphism $V\rightarrow U$ in $(\sch/S)$,
${\mathcal D}_U\ni E \mapsto E_V\in{\mathcal D}_V$
is an ``exact functor", that is, for any distinguished
triangle $E\rightarrow F \rightarrow G$ in ${\mathcal D}_U$,
$E_V\rightarrow F_V \rightarrow G_V$ is a distinguished triangle
in ${\mathcal D}_V$ and for any $E\in{\mathcal D}_U$ and any $i\in\mathbf{Z}$,
there is an isomorphism $(E[i])_V\cong E_V[i]$
functorial in $E$.
\end{enumerate}
\end{definition}

\begin{definition}\label{bcp}\rm
A fibered triangulated category $p:{\mathcal D}\rightarrow (\sch/S)$
has base change property if
\begin{enumerate}
\item for each $U\in(\sch/S)$, there is a bi-exact bi-functor
 $\otimes:{\mathcal D}_U\times D^b(\mathrm{Coh}(U/U))
 \rightarrow {\mathcal D}_U$ such that there is a functorial isomorphism 
 $E[i]\otimes P[j]\cong (E\otimes P)[i+j]$ for $E\in{\mathcal D}_U$, $P\in D^b(\mathrm{Coh}(U/U))$,
\item for a morphism $\varphi:U\rightarrow V$ in $(\sch/S)$, the diagram
\[
 \begin{CD}
  {\mathcal D}_V\times D^b(\mathrm{Coh}(V/V)) @>\otimes>> {\mathcal D}_V \\
  @V\varphi^*\times L\varphi^*VV  @VV\varphi^*V \\
  {\mathcal D}_U\times D^b(\mathrm{Coh}(U/U)) @>\otimes>> {\mathcal D}_U
 \end{CD}
\]
 is ``commutative'',
 precisely, there exists a functorial isomorphism
 $\varphi^*\circ\otimes\stackrel{\sim}\rightarrow
 \otimes\circ(\varphi^*\times L\varphi^*)$,
\item for $U\in(\sch/S)$,
 there is a bi-exact bi-functor
\[
 {\bf R}\Hom_p:{\mathcal D}_U\times{\mathcal D}_U
 \longrightarrow D^+_c(U)
\]
 such that for $E_1,E_2\in{\mathcal D}_U$ and for intgers $i,j$,
 there is an isomorphism
 $\mathbf{R}\Hom_p(E_1[i],E_2[j])\cong\mathbf{R}\Hom_p(E,F)[j-i]$
 functorial in $E_1$ and $E_2$
 and also for $E_1,E_2\in{\mathcal D}_U$ there is an isomorphism
 $\Hom_{D(U)}({\mathcal O}_U,
 {\bf R}\Hom_p(E_1,E_2))\stackrel{\sim}\rightarrow\Hom_{{\mathcal D}_U}(E_1,E_2)$
 functorial in $E_1$ and $E_2$,
\item for any $U\in(\sch/S)$ and
 for any objects $E_1,E_2\in {\mathcal D}_U$,
 there  exist a lower bounded complex $P^{\bullet}$ of locally free sheaves of finite rank
 on $U$ and an isomorphism
\[
 P^{\bullet}\otimes{\mathcal O}_V\stackrel{\sim}\rightarrow
 {\bf R}\Hom_p((E_1)_V,(E_2)_V) 
\]
in $D^+_c(V)$ for any morphism $V\rightarrow U$ in $(\sch/S)$, such that the diagram
\[
 \begin{CD}
  H^0(\Gamma((U,P^{\bullet})) & @>>> & \Hom_{D(U)}({\mathcal O}_U,\mathbf{R}\Hom_p(E_1,E_2))
  & @>\sim>> &  \Hom_{{\mathcal D}_U}(E_1,E_2) \\
  @VVV & & & &  & & @VVV \\
   H^0(\Gamma((V,P^{\bullet})) & @>>> & \Hom_{D(U)}({\mathcal O}_V,\mathbf{R}\Hom_p((E_1)_V,(E_2)_V))
   & @>\sim>> & \Hom_{{\mathcal D}_V}((E_1)_V,(E_2)_V)
 \end{CD}
\]
is commutative, 
\item  for $U\in(\sch/S)$, $E_1,E_2\in{\mathcal D}_U$
 and $F_1,F_2\in D^b(\mathrm{Coh}(U/U))$,
 there is a functorial isomorphism
 ${\bf R}\Hom_p(E_1\otimes F_1,E_2\otimes_UF_2)
 \cong{\bf R}\Hom_p(E_1,E_2)
 \otimes^L_{{\mathcal O}_U}\mathbf{R}{\mathcal Hom}(F_1,F_2)$
 such that for any morphism $\varphi:V\rightarrow U$ in $(\sch/S)$,
 the diagram
 \[
  \begin{CD}
   \mathbf{R}\Hom_p(E_1\otimes F_1, E_2\otimes F_2) @>\sim>>
   \mathbf{R}\Hom_p(E_1,E_2)\otimes^L_{{\mathcal O}_U}\mathbf{R}{\mathcal Hom}(F_1,F_2) \\
   @VVV @VVV \\
   \mathbf{R}\varphi_*(\mathbf{R}\Hom_p((E_1\otimes F_1)_V, (E_2\otimes F_2)_V)) @>\sim>>
   \mathbf{R}\varphi_*(\mathbf{R}\Hom_p((E_1)_V,(E_2)_V)\otimes^L_{{\mathcal O}_V}\mathbf{R}{\mathcal Hom}((F_1)_V,(F_2)_V)) \\
  \end{CD}
 \]
 is commutative.
\end{enumerate}
\end{definition}

\begin{remark}\rm
 For $E,F\in{\mathcal D}_U$, we denote the $i$-th cohomology
 $H^i({\bf R}\Hom_p(E,F))$ by $R^i\Hom_p(E,F)$.
 We notice that for three objects $E,F,G\in{\mathcal D}_U$,
 there is a canonical morphism
\[
 R^0\Hom_p(E,F)\times R^0\Hom_p(F,G)\rightarrow
 R^0\Hom_p(E,G).
\]
\end{remark}

\begin{example}\rm
Let $X\rightarrow S$ be a flat projective morphism.
Then $\{D^b(\mathrm{Coh}(X_U/U))\}_{U\in(\sch/S)}$
becomes a fibered triangulated category over $S$
which has base change property.
\end{example}

\begin{example}\rm
Let $X$ be a projective scheme over $\mathbf{C}$
and $G$ a finite group acting on $X$.
For a scheme $U\in(\sch/\mathbf{C})$,
let $D^G(\mathrm{Coh}(X_U/U))$ be the derived category of
bounded complexes of $G$-equivariant coherent
sheaves on $X_U$ of finite Tor-dimension over $U$.
Then $\{D^G(\mathrm{Coh}(X_U/U))\}_{U\in(\sch/\mathbf{C})}$
becomes a fibered triangulated category over $\mathbf{C}$
which has base change property.
\end{example}

\section{Strict ample sequence and stability}

\begin{definition}\label{def-polarization}\rm
Let $p:{\mathcal D}\rightarrow(\sch/S)$
be a fibered triangulated category with base change property.
A sequence ${\mathcal L}=\{L_n\}_{n\geq 0}$ of objects of
${\mathcal D}_S$ is said to be a strict ample sequence
if it satisfies the following conditions:
\begin{enumerate}
\item $\Ext^i((L_N)_s,(L_n)_s)=0$
 for any $i\neq 0$, $N>n$ and $s\in S$.
\item There exist isomorphisms
\[
 \theta_k: R^0\Hom_p(L_n,L_m)
 \stackrel{\sim}\longrightarrow
 R^0\Hom_p(L_{n+k},L_{m+k})
\]
 for non-negative integers $k,m,n$ with $n\geq m$ such that
 $\theta_k\circ\theta_l=\theta_{k+l}$ for any $k,l$ and the diagram
\[
\begin{CD}
 R^0\Hom_p(L_n,L_m)\otimes R^0\Hom_p(L_m,L_l)
 @>\theta_k\otimes\theta_k>>
 R^0\Hom_p(L_{n+k},L_{m+k})\otimes R^0\Hom_p(L_{m+k},L_{l+k}) \\
 @VVV @VVV \\
 R^0\Hom_p(L_n,L_l) @>\theta_k>> R^0\Hom_p(L_{n+k},L_{l+k})
\end{CD}
\]
 is commutative for non-negative integers $k,l,m,n$ with
 $n\geq m\geq l$.
\item There exists a subbundle $V_1\subset R^0\Hom_p(L_1,L_0)$
such that the diagram
\[
 \begin{CD}
  V_1\times R^0\Hom_p(L_n,L_0)
  @>\theta_n\times 1>>
  R^0\Hom_p(L_{n+1},L_n)\times R^0\Hom_p(L_n,L_0) \\
  @V 1\times\theta_1 VV  @VVV \\
  V_1\times R^0\Hom_p(L_{n+1},L_1)
  @>>> R^0\Hom_p(L_{n+1},L_0),
 \end{CD}
\]
is commutative for $n\geq 0$,
where the right vertical arrow and the bottom
horizontal arrow are the canonical composition maps
and there exists an integer $n_0$ such that for any $n\geq n_0$,
 \[
  R^0\Hom_p(L_n,L_1)\otimes V_1 \longrightarrow
  R^0\Hom_p(L_n,L_0)
 \]
is surjective for any $n\geq n_0$.
\item For any object $E\in{\mathcal D}_U$ and for any non-negative integer $m$,
there exists a bounded complex $P^{\bullet}$ of locally free sheaves of finite rank on $U$ such that
$\mathbf{R}\Hom_p((L_m)_V,E_V)\cong P^{\bullet}\otimes{\mathcal O}_V$ for any $V\rightarrow U$. 
Moreover, there exists an integer $n_0$ such that for any $n\geq n_0$,
 exists an integer $N_0$ such that
 for any integers $i,N$ with $N\geq N_0$ and for any $s\in U$,
 \[
  \Hom((L_N)_s,(L_n)_s)\otimes \Ext^i((L_n)_s,E_s)\rightarrow
  \Ext^i((L_N)_s,E_s)
 \]
 is surjective.
\item If there exist integers $i$, $n_0$ and
 an object $E\in{\mathcal D}_U$ satisfying 
 $\Ext^i((L_n)_s,E_s)=0$ for any $n\geq n_0$ and for any $s\in U$,
 then there exist an object $F\in {\mathcal D}_U$ and a morphism
 $u:E\rightarrow F$ such that for any $j>i$,
 $R^j\Hom_p((L_n)_U,E)\rightarrow
 R^j\Hom_p((L_n)_U,F)$
 are isomorphic for $n\gg 0$,
 and for any $j\leq i$, $R^j\Hom_p((L_n)_U,F)=0$ for $n\gg 0$.
\item Take two objects $E,F\in{\mathcal D}_U$ such that for any $i\geq 0$,
 $R^i\Hom_p((L_n)_U,E)=0$ for $n\gg 0$
 and that for any $i<0$,
 $R^i\Hom_p((L_n)_U,F)=0$ for $n\gg 0$.
 Then we have
 $\Hom_{{\mathcal D}_U}(E,F)=0$.
\end{enumerate}
\end{definition}

\begin{proposition}\label{prop:object-vanish}
 Take $E\in{\mathcal D}_U$ such that for any $i$,
 $R^i\Hom_p((L_n)_U,E)=0$ for $n\gg 0$.
 Then we have $E=0$.
\end{proposition}

\begin{proof}
Applying the condition (6) of Definition \ref{def-polarization},
we have $\Hom(E,E)=0$.
In particular $\mathrm{id}_E=0$.
So, for any object $F\in{\mathcal D}_U$ and for any morphism
$f\in\Hom(F,E)$ (resp.\ $g\in\Hom(E,F)$),
$f=\mathrm{id}_E\circ f=0$ (resp.\ $g=g\circ\mathrm{id}_E=0$).
Thus $E=0$.
\end{proof}

\begin{remark}\rm
By the condition in Definition \ref{def-polarization} (2),
we can see that $\theta_0=\mathrm{id}$ and
$\theta_k(\mathrm{id})=\mathrm{id}$.
We put $A:=\bigoplus_{n\geq 0} R^0\Hom_p(L_n,L_0)$
and define a multiplication
\[
 \alpha:R^0\Hom_p(L_n,L_0)\times R^0\Hom_p(L_m,L_0)
 \longrightarrow R^0\Hom_p(L_{n+m},L_0)
\]
by $\alpha=(\text{composition})\circ(\theta_m\times\mathrm{id})$.
Then $A$ becomes an associative graded ring
which is a finitely generated module over $S^*(V_1)$,
where $S^*(V_1)$ is the symmetric algebra of $V_1$ over ${\mathcal O}_S$.
\end{remark}

\begin{proposition}\label{isom-n>>0}
Let $E_1,E_2$ be objects of ${\mathcal D}_U$ and
$u:E_1\rightarrow E_2$ be a morphism such that for any integer $i$
the induced morphism
$R^i\Hom_p((L_n)_U,E_1)\rightarrow
R^i\Hom_p((L_n)_U,E_2)$
is isomorphic for $n\gg 0$.
Then $u$ is an isomorphism.
\end{proposition}

\begin{proof}
For any $i$, there is an exact sequence
\begin{gather*}
 R^i\Hom_p((L_n)_U,E_1)\stackrel{\sim}\longrightarrow
 R^i\Hom_p((L_n)_U,E_2) \longrightarrow
 R^i\Hom_p((L_n)_U,\mathrm{Cone}(u)) \\
 \longrightarrow R^{i+1}\Hom_p((L_n)_U,E_1)\stackrel{\sim}\longrightarrow
 R^{i+1}\Hom_p((L_n)_U,E_2)
\end{gather*}
for $n\gg 0$.
Thus we have $R^i\Hom_p((L_n)_U,\mathrm{Cone}(u))=0$ for $n\gg 0$.
By Proposition \ref{prop:object-vanish} we have
$\mathrm{Cone}(u)=0$, which means that $u$ is an isomorphism.
\end{proof}

\begin{proposition}
For an integer $i$ and an object $E\in{\mathcal D}_U$ such that for $n\gg 0$,
$\Ext^i((L_n)_s,E_s)=0$ for $s\in U$, the object $F$ given in
Definition \ref{def-polarization} (5) is unique up to an isomorphism.
\end{proposition}

\begin{proof}
Let $F'\in{\mathcal D}_U$ be another object with a morphism
$u':E\rightarrow F'$ having the same property as $F$.
Consider the composite
\[
 v:\mathrm{Cone}(u)[-1]\longrightarrow E \stackrel{u'}\longrightarrow F'.
\]
Since there is a long exact sequence
\begin{gather*}
 \cdots\longrightarrow R^j\Hom_p((L_n)_U,E) \longrightarrow R^j\Hom_p((L_n)_U,F) 
 \longrightarrow R^j\Hom_p((L_n)_U,\mathrm{Cone}(u)) \\
 \longrightarrow R^{j+1}\Hom_p((L_n)_U,E)\longrightarrow
 R^{j+1}\Hom_p((L_n)_U,F)\longrightarrow \cdots, 
\end{gather*}
we have, for any $j\geq i$,
$R^j\Hom_p((L_n)_U,\mathrm{Cone}(u))=0$ for $n\gg 0$.
Note that for any $j\leq i$, we have $R^j\Hom_p((L_n)_U,F')=0$ for $n\gg 0$.
Then we have $\Hom_{{\mathcal D}_U}(\mathrm{Cone}(u),F')=0$ and
$\Hom_{{\mathcal D}_U}(\mathrm{Cone}(u)[-1],F')=0$
by condition (6) of Definition \ref{def-polarization}.
So we have $v=0$ and there is a unique morphism $\varphi:F\rightarrow F'$
which makes the diagram
\[
 \begin{CD}
  E @>u>> F \\
  @V\mathrm{id}VV @VV\varphi V \\
  E @>u'>> F'
 \end{CD}
\]
commute.
We can see that for any integer $j$, the morphism
$R^j\Hom_p((L_n)_U,F)\rightarrow R^j\Hom_p((L_n)_U,F')$
induced by $\varphi$ is isomorphic for $n\gg 0$.
Hence $\varphi$ is an isomorphism by Proposition \ref{isom-n>>0}
\end{proof}

\begin{remark}\label{flattening}\rm
In the situation of Definition \ref{def-polarization} (5),
for $n\gg 0$, the induced morphism
\[
 \Ext^j((L_n)_s,E_s)\rightarrow\Ext^j((L_n)_s,F_s)
\]
is isomorphic for any $j>i$ and for any $s\in U$,
and we have, for $n\gg 0$, 
$\Ext^j((L_n)_s,F_s)=0$ for any $j\leq i$ and for any $s\in U$.

Indeed consider the distinguished triangle
$E\stackrel{u}\rightarrow F \rightarrow \mathrm{Cone}(u)$.
Note that there is a long exact sequence
\begin{gather*}
 R^j\Hom_p((L_n)_U,E)\longrightarrow R^j\Hom_p((L_n)_U,F)\longrightarrow
 R^j\Hom_p((L_n)_U,\mathrm{Cone}(u)) \\
 \longrightarrow R^{j+1}\Hom_p((L_n)_U,E)\longrightarrow
 R^{j+1}\Hom_p((L_n)_U,F).
\end{gather*}
Since $R^i\Hom_p((L_n)_U,F)=0$ for $n\gg 0$,
and for any $j>i$, $R^j\Hom_p((L_n)_U,E)\rightarrow R^j\Hom_p((L_n)_U,F)$
are isomorphic for $n\gg 0$, we have, for any $j\geq i$,
$R^j\Hom_p((L_n)_U,\mathrm{Cone}(u))=0$ for $n\gg 0$.

By Definition \ref{def-polarization} (4),
there are integers $n_0$ and $N_0$ with $N_0>n_0$ such that
\[
 \Hom((L_N)_s,(L_{n_0})_s)\otimes\Ext^j((L_{n_0})_s,\mathrm{Cone}(u)_s)
 \longrightarrow \Ext^j((L_N)_s,\mathrm{Cone}(u)_s)
\]
is suriective for any $j$, any $N\geq N_0$ and any $s\in U$.
By Definition \ref{def-polarization} (4), there are integers $j_0,j_1$ such that
for $j<j_0$ and $j>j_1$, $\Ext^j((L_{n_0})_s,\mathrm{Cone}(u)_s)=0$ for any $s\in U$.
Then for any $N\geq N_0$, we have $\Ext^j((L_N)_s,\mathrm{Cone}(u)_s)=0$
for any $j>j_1$ and $s\in U$.
For each $j$ with $i\leq j\leq j_1$, there exists an integer $N(j)$ such that
for any $N\geq N(j)$, we have $R^j\Hom_p((L_N)_U,\mathrm{Cone}(u))=0$.
Put 
\[
 \tilde{N}:=\max\{N(i),N(i+1),\ldots,N(j_1),N_0\}.
 \]
By Definition \ref{bcp} (4),
we have $\Ext^j((L_N)_s,\mathrm{Cone}(u)_s)=0$
for any $N\geq \tilde{N}$ and for each $j$ with $i\leq j\leq j_1$ and for any $s\in U$,
because $\Ext^{j_1+1}((L_N)_s,\mathrm{Cone}(u)_s)=0$ for any $s\in U$ and
$R^j\Hom_p((L_N)_U,\mathrm{Cone}(u))=0$ for $i\leq j\leq j_1$.
Thus we have $\Ext^j((L_N)_s,\mathrm{Cone}(u)_s)=0$
for any $N\geq\tilde{N}$, $j\geq i$ and $s\in U$.

Note that there are integers $k_0,k_1$ and a positive integer $M_0$ such that
for any $M\geq M_0$ and for any $s\in U$,
$\Ext^j((L_M)_s,F_s)=0$ for $j<k_0$ and $j>k_1$.
We may also assume that for any $M\geq M_0$ and for any $s\in U$,
$\Ext^i((L_M)_s,E_s)=0$.
From the exact sequence
\[
 0=\Ext^i((L_M)_s,E_s) \longrightarrow \Ext^i((L_M)_s,F_s)
 \longrightarrow \Ext^i((L_M)_s,\mathrm{Cone}(u)_s)=0,
\]
we have $\Ext^i((L_M)_s,F_s)=0$ for $s\in U$ and $M\geq \max\{M_0,\tilde{N}\}$.
By assumption, for each $j$ with $k_0\leq j\leq i$, there exists an integer $M(j)$ such that
$R^j\Hom_p((L_M)_U,F)=0$ for $M\geq M(j)$.
Put
\[
 \tilde{M}:=\max\{\tilde{N},M_0,M(k_0),M(k_0+1),\ldots,M(i)\}.
 \]
Then we have $\Ext^j((L_M)_s,F_s)=0$ for $j\leq i$, $s\in U$ and $M\geq \tilde{M}$
by using Definition \ref{bcp} (4),
because $\Ext^i((L_M)_s,F_s)=0$ and $R^j\Hom_p((L_M)_U,F)=0$
for $k_0\leq j\leq i$.
From the exact sequence
\[
 \Ext^{j-1}((L_M)_s,\mathrm{Cone}(u)_s) \longrightarrow
 \Ext^j((L_M)_s,E_s) \longrightarrow \Ext^j((L_M)_s,F_s)
 \longrightarrow \Ext^j((L_M)_s,\mathrm{Cone}(u)_s),
\]
we have an isomorphism
$\Ext^j((L_M)_s,E_s) \stackrel{\sim}\rightarrow \Ext^j((L_M)_s,F_s)$
for $j>i$, $s\in U$ and $M\geq \tilde{M}$.
\end{remark}

\begin{lemma}\label{lemma:resolution1}
 If $E\in{\mathcal D}_U$ satisfies $\Ext^i((L_n)_s,E_s)=0$
 for $n\gg 0$, $i\neq 0$ and $s\in U$, then
 there exist locally free ${\mathcal O}_U$-modules
 $W_0,W_1,W_2$, positive integers $n_0<n_1<n_2$ and
 morphisms
\[
 (L_{n_2})_U\otimes W_2\stackrel{d^2}\longrightarrow
 (L_{n_1})_U\otimes W_1\stackrel{d^1}\longrightarrow
 (L_{n_0})_U\otimes W_0\stackrel{f}\longrightarrow E
\]
 such that the induced sequence
\begin{gather*}
 \Hom((L_N)_s,(L_{n_2})_s)\otimes W_2\longrightarrow
 \Hom((L_N)_s,(L_{n_1})_s)\otimes W_1\longrightarrow \\
 \Hom((L_N)_s,(L_{n_0})_s)\otimes W_0\longrightarrow
 \Hom((L_N)_s,E_s)\longrightarrow 0
\end{gather*}
is exact for $N\gg 0$ and $s\in U$.
\end{lemma}

\begin{proof}
By Definition \ref{def-polarization} (4), there exist integers
$n_0,N_0$ with $N_0>n_0$ such that for any $s\in U$,
\[
 \Hom((L_N)_s,(L_{n_0})_s)\otimes\Hom((L_{n_0})_s,E_s)\rightarrow
 \Hom((L_N)_s,E_s)
\]
is surjective for $N\geq N_0$ and
$\Ext^i((L_n)_s,E_s)=0$ for $n\geq n_0$, $i\neq 0$ and $s\in U$.
There is a canonical morphism
\[
 f:(L_{n_0})_U\otimes R^0\Hom_p((L_{n_0})_U,E) \longrightarrow E
\]
and we put $F^1:=\mathrm{Cone}(f)[-1]$.
Then we can see that $\Ext^i((L_N)_s,(F^1)_s)=0$ for $N\geq N_0$, $i\neq 0$
and $s\in U$.
We can find integers $n_1,N_1$ with $N_1>n_1$ such that for any $s\in U$,
\[
 \Hom((L_N)_s,(L_{n_1})_s)\otimes\Hom((L_{n_1})_s,(F^1)_s)
 \longrightarrow \Hom((L_N)_s,(F^1)_s)
\]
is surjective for $N\geq N_1$ and
$\Ext^i((L_n)_s,(F^1)_s)=0$ for $n\geq n_1$, $i\neq 0$ and $s\in U$.
Consider the canonical morphism
\[
 g:(L_{n_1})_U\otimes R^0\Hom_p((L_{n_1})_U,F^1)
 \longrightarrow F^1
\]
and put $F^2:=\mathrm{Cone}(g)[-1]$.
We can find again integers $n_2,N_2$ with $N_2>n_2$ such that for any $s\in U$,
\[
 \Hom((L_N)_s,(L_{n_2})_s)\otimes\Hom((L_{n_2})_s,(F^2)_s)
 \longrightarrow \Hom((L_N)_s,(F^2)_s)
\]
is surjective for $N\geq N_2$ and
$\Ext^i((L_n)_s,(F^2)_s)=0$ for $n\geq n_2$, $i\neq 0$ and $s\in U$.
There is a canonical morphism
\[
 h:(L_{n_2})_U\otimes R^0\Hom_p((L_{n_2})_U,F^2)
 \longrightarrow F^2
\]
and we obtain a sequence of morphisms
\begin{gather*}
 (L_{n_2})_U\otimes R^0\Hom_p((L_{n_2})_U,F^2) \longrightarrow
 (L_{n_1})_U\otimes R^0\Hom_p((L_{n_1})_U,F^1) \\
 \longrightarrow
 (L_{n_0})_U\otimes R^0\Hom_p((L_{n_0})_U,E) \longrightarrow E
\end{gather*}
such that for $N\geq \max\{N_0,N_1,N_2\}$, the induced sequence
\begin{gather*}
 \Hom((L_N)_s,(L_{n_2})_s)\otimes R^0\Hom_p((L_{n_2})_U,F^2)\longrightarrow
 \Hom((L_N)_s,(L_{n_1})_s)\otimes R^0\Hom_p((L_{n_1})_U,F^1) \\
 \longrightarrow \Hom((L_N)_s,(L_{n_0})_s)\otimes R^0\Hom_p((L_{n_0})_U,E)
 \longrightarrow \Hom((L_N)_s,E_s) \longrightarrow 0
\end{gather*}
is exact for any $s\in U$.
If we put $W_0=R^0\Hom_p((L_{n_0})_U,E)$ and $W_i=R^0\Hom_p((L_{n_i})_U,F^i)$ for $i=1,2$,
then we can see by Definition \ref{bcp} (4) that $W_i$ are locally free ${\mathcal O}_U$-modules
and have the desired property.
\end{proof}

\begin{proposition}\label{prop:fully-faithful}
Let $E_1,E_2$ be objects of ${\mathcal D}_U$ such that
$\Ext^i((L_n)_s,(E_j)_s)=0$ for $j=1,2$, $n\gg 0$, $i\neq 0$ and $s\in U$.
If $f:E_1\rightarrow E_2$ is a morphism in ${\mathcal D}_U$
such that the induced morphisms
$R^0\Hom_p((L_n)_U,E_1)\rightarrow R^0\Hom_p((L_n)_U,E_2)$
are zero for $n\gg0$, then $f=0$.
\end{proposition}

\begin{proof}
By assumption, there is an integer $N_0$ such that for any $N\geq N_0$,
the morphism
\[
 R^0\Hom_p((L_N)_U,E_1)\rightarrow R^0\Hom_p((L_N)_U,E_2)
\]
induced by $f$ is zero
and $\Ext^i((L_N)_s,(E_j)_s)=0$ for $j=1,2$, $i\neq 0$ and $s\in U$.
By Lemma \ref{lemma:resolution1}, there are locally free sheaves
$W_0,W_1,W_2$, integers $n_0<n_1<n_2$ and morphisms
\[
 (L_{n_2})_U\otimes W_2 \longrightarrow (L_{n_1)_U}\otimes W_1
 \longrightarrow (L_{n_0)_U}\otimes W_0 \stackrel{\varphi}\longrightarrow E_1
\]
such that the induced sequence
\begin{gather*}
 \Hom((L_N)_s,(L_{n_2})_s)\otimes W_2\longrightarrow
 \Hom((L_N)_s,(L_{n_1})_s)\otimes W_1\longrightarrow \\
 \Hom((L_N)_s,(L_{n_0})_s)\otimes W_0\longrightarrow
 \Hom((L_N)_s,(E_1)_s)\longrightarrow 0
\end{gather*}
is exact for $N\gg 0$ and $s\in U$.
We can take $n_0$ so that $n_0\geq N_0$.
Consider the distinguished triangle
\[
 (L_{n_0})_U\otimes W_0 \longrightarrow E_1
 \longrightarrow \mathrm{Cone}(\varphi).
\]
We can see that
$\Ext^i((L_n)_s,\mathrm{Cone}(\varphi)_s)=0$ for $n\gg 0$, $i\neq -1$
and $s\in U$.
So we have $\Hom_{{\mathcal D}_U}(\mathrm{Cone}(\varphi),E_2)=0$
by (6) of Definition \ref{def-polarization} and the homomorphism
\[
 (\dag) \quad
 \Hom_{{\mathcal D}_U}(E_1,E_2)\rightarrow
 \Hom_{{\mathcal D}_U}((L_{n_0})_U\otimes W_0,E_2)
\]
induced by $\varphi$ is injective.
On the other hand, the homomorphism
\[
 R^0\Hom_p((L_{n_0})_U\otimes W_0,E_1) \longrightarrow 
 R^0\Hom_p((L_{n_0})_U\otimes W_0,E_2)
\]
induced by $f$ is zero.
So we have $f\circ\varphi=0$.
By the injectivity of $(\dag)$, we have $f=0$.
\end{proof}

Since $A=\bigoplus_{n\geq 0}R^0\Hom_p(L_n,L_0)$ becomes
a finite algebra over $S^*(V_1)$,
the associated sheaf ${\mathcal A}:=\tilde{A}$ becomes
a coherent sheaf of algebras on $\mathbf{P}(V_1)$.
For each object $E\in{\mathcal D}_U$ 
satisfying $\Ext^i((L_n)_s,E_s)=0$ for $n\gg 0$, $i\neq 0$ and $s\in U$,
the associated sheaf
$\left(\bigoplus_{n\geq 0} R^0\Hom_p((L_n)_U,E)\right)^{\sim}$
on $\mathbf{P}(V_1)_U=\mathbf{P}(V_1)\times_SU$
becomes a coherent ${\mathcal A}_U$-module flat over $U$.

\begin{proposition}\label{equivalence}
 The correspondence
 $E\mapsto \left(\bigoplus_{n\geq 0} R^0\Hom_p((L_n)_U,E)\right)^{\sim}$
 gives an equivalence of categories between the full subcategory of
 ${\mathcal D}_U$ consisting of the objects $E$ of ${\mathcal D}_U$
 satisfying $\Ext^i((L_n)_s,E_s)=0$ for $n\gg 0$, $i\neq 0$ and $s\in U$
 and the category of coherent ${\mathcal A}_U$-modules flat over $U$.
\end{proposition}

\begin{proof}
First we will prove that the functor
\[
 \psi: E \mapsto \left( \bigoplus_{n\geq 0} R^0\Hom_p((L_n)_U,E)\right)^{\sim}
\]
is fully faithful.
Take any objects $E,F$ of ${\mathcal D}_U$ which satisfy $\Ext^i((L_n)_s,E_s)=0$,  $\Ext^i((L_n)_s,F_s)=0$
for $n\gg 0$, $i\neq 0$ and $s\in U$. 
By Proposition \ref{prop:fully-faithful},  
\[
 (\dag) \quad \Hom(E,F)\longrightarrow \Hom(\psi(E),\psi(F))
\]
is injective.
Take any homomorphism
$f \in \Hom(\psi(E),\psi(F))$.
There exists an integer $n_0$ such that for any $n\geq n_0$,
$\Ext^i((L_n)_s,E_s)=0$, $\Ext^i((L_n)_s,F_s)=0$ for $i\neq 0$ and $s\in U$
and the homomorphisms
\begin{gather*}
 \Hom((L_N)_s,(L_{n_0})_s)\otimes\Hom((L_{n_0})_s,E_s)\longrightarrow\Hom((L_N)_s,E_s) \\
 \Hom((L_N)_s,(L_{n_0})_s)\otimes\Hom((L_{n_0})_s,F_s)\longrightarrow\Hom((L_N)_s,F_s) 
\end{gather*}
are surjective for $N\gg 0$ and $s\in U$.
For a coherent ${\mathcal A}_U$-module ${\mathcal E}$, we denote
${\mathcal E}\otimes{\mathcal O}_{\mathbf{P}(V_1)_U}(n)$ simply by ${\mathcal E}(n)$.
We denote the structure morphism $\mathbf{P}(V_1)_U\rightarrow U$ by $\pi$.
Then we may assume that
$R^i\pi_*(\psi(E)(n_0))=0$, $R^i\pi_*(\psi(F)(n_0))=0$ for $i>0$
and that the homomorphisms
\begin{gather*}
 \pi_*(\psi(E)(n_0))\otimes{\mathcal A}(-n_0)\longrightarrow \psi(E) \\
 \pi_*(\psi(F)(n_0))\otimes{\mathcal A}(-n_0)\longrightarrow \psi(F)
\end{gather*}
are surjective.
We may also assume that
\begin{gather*}
 R^0\Hom((L_{n_0})_U,E) \longrightarrow \pi_*(\psi(E)(n_0))  \\
 R^0\Hom((L_{n_0})_U,F) \longrightarrow \pi_*(\psi(F)(n_0))
 \end{gather*}
 are isomorphic.
Consider the distinguished triangles
\begin{gather*}
 \mathrm{Cone}(v)[-1] \stackrel{\iota_1}\longrightarrow (L_{n_0})_U\otimes R^0\Hom_p((L_{n_0})_U,E)
 \stackrel{v}\longrightarrow E   \\
 \mathrm{Cone}(w)[-1] \stackrel{\iota_2}\longrightarrow (L_{n_0})_U\otimes R^0\Hom_p((L_{n_0})_U,F)
 \stackrel{w}\longrightarrow F.
\end{gather*}
Then we can see that $\Ext^i((L_N)_s,\mathrm{Cone}(v)_s)=0$, $\Ext^i((L_N)_s,\mathrm{Cone}(w)_s)=0$
for $N\gg 0$, $i\neq 0$ and $s\in U$.
The homomorphism $f:\psi(E)\rightarrow\psi(F)$ induces a homomorphism
\[
 f(n_0):R^0\Hom_p((L_{n_0})_U,E) \cong \pi_*(\psi(E)(n_0)) \longrightarrow
 \pi_*(\psi(F)(n_0)) \cong R^0\Hom((L_{n_0})_U,F).
\]
Then $f(n_0)$ induces a homomorphism
\[
 \tilde{f}:(L_{n_0})_U\otimes R^0\Hom_p((L_{n_0})_U,E)\longrightarrow (L_{n_0})_U\otimes R^0\Hom_p((L_{n_0})_U,F).
\]
Consider the composite
\[
 w\circ\tilde{f}\circ\iota_1: \mathrm{Cone}(v)[-1] \stackrel{\iota_1}\longrightarrow (L_{n_0})_U\otimes \Hom((L_{n_0})_U,E)
 \stackrel{\tilde{f}}\longrightarrow (L_{n_0})_U\otimes \Hom((L_{n_0})_U,F)
 \stackrel{w}\longrightarrow F.
\]
Then we have
$\psi(w\circ\tilde{f}\circ\iota_1)=\psi(w)\circ\psi(\tilde{f})\circ\psi(\iota_1)=f\circ \psi(v)\circ\psi(\iota_1)
=f\circ\psi(v\circ\iota_1)=0$.
Since
\[
 \Hom(\mathrm{Cone}(v)[-1],F)\longrightarrow \Hom(\psi(\mathrm{Cone}(v)[-1]),\psi(F))
\]
is injective, we have $w\circ\tilde{f}\circ\iota_1=0$.
So there is a morphism $f':E\rightarrow F$, which makes the diagram
\[
\begin{CD}
 \mathrm{Cone}(v)[-1] & @>\iota_1>> & (L_{n_0})_U\otimes R^0\Hom_p((L_{n_0})_U,E) & @>v>> E \\
  & & & & @V\tilde{f}VV & @V f' VV \\
 \mathrm{Cone}(w)[-1] & @>\iota_2>> & (L_{n_0})_U\otimes R^0\Hom_p((L_{n_0})_U,F) & @>w>> F
\end{CD}
\]
commute.
This commutative diagram induces a commutative diagram
\[
\begin{CD}
  \pi_*(\psi(E)(n_0))\otimes{\mathcal A}(-n_0) @>\psi(v)>> \psi(E) \\
  @VVV @V\psi(f') VV \\
  \pi_*(\psi(F)(n_0))\otimes{\mathcal A}(-n_0) @>>> \psi(F).
\end{CD}
\]
Since
$\psi(v)\circ(\psi(f')-f)=\psi(v)\circ\psi(f')-\psi(v)\circ f
=\psi(w)\circ\psi(\tilde{f})-\psi(w)\circ\psi(\tilde{f})=0$,
we have $\psi(f')-f=0$ because $\psi(v)$ is surjective.
So we have $\psi(f')=f$.
Thus $(\dag)$ is surjective and $\psi$ becomes a fully faithful functor.

Take any coherent ${\mathcal A}_U$-module ${\mathcal E}$ flat over $U$.
There is an exact sequence of coherent ${\mathcal A}_U$-modules
\[
 W_2\otimes {\mathcal A}(-n_2) \stackrel{\delta^2}\longrightarrow W_1\otimes{\mathcal A}(-n_1)
 \stackrel{\delta^1}\longrightarrow W_0\otimes{\mathcal A}(-n_0) \longrightarrow {\mathcal E}\longrightarrow 0,
\]
where $W_0, W_1, W_2$ are locally free sheaves on $U$ and $n_2\gg n_1\gg n_0\gg 0$.
The above sequence induces a sequence of morphisms
\[
 (L_{n_2})_U\otimes W_2 \stackrel{d^2}\longrightarrow (L_{n_1})_U\otimes W_1
 \stackrel{d^1}\longrightarrow (L_{n_0})_U\otimes W_0.
\]
By construction we have $d^1\circ d^2=0$.
So there is a morphism $u:\mathrm{Cone}(d^2)\rightarrow (L_{n_0})_U\otimes W_0$
such that the diagram
\[
 \begin{array}{ccc}
   (L_{n_1})_U\otimes W_1 & \stackrel{d^1}\longrightarrow &  (L_{n_0})_U\otimes W_0 \\
   \searrow \hspace{-50pt} & & \hspace{-40pt} \nearrow _u \\
   & \mathrm{Cone}(d^2) &
 \end{array}
\]
is commutative.
Note that $\Ext^i((L_N)_s,\mathrm{Cone}(d^2)_s)=0$ for $N\gg 0$, $i\neq -1,0$ and $s\in U$.
So we have $\Ext^i((L_N)_s,\mathrm{Cone}(u)_s)=0$ for $N\gg 0$, $i\neq -2,-1,0$ and $s\in U$.
Since ${\mathcal E}$ is flat over $U$, the sequence
\[
 W_2\otimes{\mathcal A}(-n_2)\otimes k(s)  \longrightarrow W_1\otimes{\mathcal A}(-n_1)\otimes k(s)
 \longrightarrow W_0\otimes{\mathcal A}(-n_0)\otimes k(s) \longrightarrow {\mathcal E}\otimes k(s) \longrightarrow 0
\]
is exact for any $s\in U$.
So we obtain the exact commutative diagram
\[
\begin{CD}
 H^0(W_2\otimes{\mathcal A}(N-n_2)\otimes k(s)) @>>> H^0(W_1\otimes{\mathcal A}(N-n_1)\otimes k(s))
 @>>> H^0(W_0\otimes{\mathcal A}(N-n_0)\otimes k(s)) \\
 @V\cong VV @V\cong VV @V\cong VV \\
 \Hom((L_N)_s,(L_{n_2})_s\otimes (W_2)_s) @>>> \Hom((L_N)_s,(L_{n_1})_s\otimes (W_1)_s) @>>>
 \Hom((L_N)_s,(L_{n_0})_s\otimes (W_0)_s)
\end{CD}
\]
for $N\gg 0$ and $s\in U$.
Here we denote $W_i\otimes k(s)$ by $(W_i)_s$ for $i=0,1,2$.
We have a factorization
\[
 \begin{array}{ccc}
  \Hom((L_N)_s,(L_{n_1})_s\otimes (W_1)_s) & \longrightarrow &  \Hom((L_N)_s,(L_{n_0})_s\otimes (W_0)_s) \\
  \hspace{80pt} \searrow & & \nearrow \hspace{80pt} \\
  & \Hom((L_N)_s,\mathrm{Cone}(d^2)_s) &
 \end{array}
\]
for $N\gg 0$ and $s\in U$,
and the homomorphism 
$\Hom((L_N)_s,(L_{n_1})_s\otimes (W_1)_s)\longrightarrow \Hom((L_N)_s,\mathrm{Cone}(d^2)_s)$
is surjective for $N\gg 0$ and $s\in U$, because
$\Ext^1((L_N)_s,(L_{n_2})_s\otimes (W_2)_s)=0$ for $N\gg 0$ and $s\in U$.
So we can see that the homomorphism
\[
  \Hom((L_N)_s,\mathrm{Cone}(d^2)_s) \longrightarrow
  \Hom((L_N)_s,(L_{n_0})_s\otimes (W_0)_s)
\]
 is injective for $N\gg 0$ and $s\in U$. 
 Since there is an exact sequence
 \begin{gather*}
  0=\Ext^{-1}((L_N)_s,(L_{n_0})_s\otimes (W_0)_s) \longrightarrow \Ext^{-1}((L_N)_s,\mathrm{Cone}(u)_s)  \\
  \stackrel{0}\longrightarrow
  \Hom((L_N)_s,\mathrm{Cone}(d^2)_s) \longrightarrow \Hom((L_N)_s,(L_{n_0})_s\otimes (W_0)_s)
 \end{gather*}
for $N\gg 0$ and $s\in U$,
we have
$\Ext^{-1}((L_N)_s,\mathrm{Cone}(u)_s)=0$ for $N\gg 0$ and $s\in U$.
By Definition \ref{def-polarization} (5) and Remark \ref{flattening},
there is an object $E\in{\mathcal D}_U$ and a morphism $\alpha:\mathrm{Cone}(u)\rightarrow E$
such that 
$R^0\Hom_p((L_N)_U,\mathrm{Cone}(u)) \rightarrow R^0\Hom_p((L_N)_U,E)$
is isomorphic for $N\gg 0$ and that
$\Ext^j((L_N)_s,E_s)=0$ for $N\gg 0$, $j\neq 0$ and $s\in U$.
We can see that the sequence
\begin{gather*}
 R^0\Hom_p((L_N)_U,(L_{n_2})_U\otimes W_2) \longrightarrow R^0\Hom_p((L_N)_U,(L_{n_1})_U\otimes W_1) \longrightarrow  \\
 R^0\Hom_p((L_N)_U,(L_{n_0})_U\otimes W_0) \longrightarrow R^0\Hom_p((L_N)_U,\mathrm{Cone}(u)) \longrightarrow 0
\end{gather*}
is exact.
Since $R^0\Hom_p((L_N)_U,\mathrm{Cone}(u))\cong R^0\Hom_p((L_N)_U,E)$ for $N\gg 0$,
there is an integer $N_0$ such that for any $N\geq N_0$,
there is a unique isomorphism $R^0\Hom_p((L_N)_U,E)\stackrel{\sim}\rightarrow \pi_*({\mathcal E}(N))$
which makes the diagram
\[
 \begin{CD}
  R^0\Hom_p((L_N)_U,(L_{n_1})_U\otimes W_1) @>>> R^0\Hom_p((L_N)_U,(L_{n_0})_U\otimes W_0) @>>>
  R^0\Hom_p((L_N)_U,E) \longrightarrow 0 \\
  @V\cong VV @V\cong VV @V\cong VV \\
  \pi_*({\mathcal A}(N-n_1)\otimes W_1) @>>>  \pi_*({\mathcal A}(N-n_0)\otimes W_0)
  @>>> \pi_*({\mathcal E}(N)) \longrightarrow 0
 \end{CD}
\]
commute.
Note that there is a canonical commutative diagram
\[
 \begin{CD}
  R^0\Hom_p((L_{N+m})_U,(L_N)_U)\otimes R^0\Hom_p((L_N)_U,E) @>>> R^0\Hom_p((L_{N+m})_U,E) \\
  @VVV @VVV \\
  \pi_*({\mathcal A}(m))\otimes \pi_*({\mathcal E}(N)) @>>> \pi_*({\mathcal E}(N+m))
 \end{CD}
\]
for $N\geq N_0$ and a non-negative integer $m$.
Then we have an isomorphism
\[
 \bigoplus_{n\geq N_0} R^0\Hom_p((L_n)_U,E) \stackrel{\sim}\longrightarrow
 \bigoplus_{n\geq N_0} \pi_*({\mathcal E}(n))
\]
of graded $A_U$-modules.
So we obtain an isomorphism
\[
 \psi(E)=\left( \bigoplus_{n\geq N_0} R^0\Hom_p((L_n)_U,E) \right)^{\sim}
 \stackrel{\sim}\longrightarrow
 \left( \bigoplus_{n\geq N_0} \pi_*({\mathcal E}(n)) \right)^{\sim}
 \cong {\mathcal E}.
\]
Thus $\psi$ becomes an equivalence of categories.
\end{proof}

\begin{definition}\rm
For a geometric point $\Spec k \rightarrow S$, an object
$E\in{\mathcal D}_k$ is said to be ${\mathcal L}$-stable
(resp.\ ${\mathcal L}$-semistable) if 
$\Ext^i((L_n)_k,E)=0$ for $n\gg 0$ and $i\neq 0$
and the inequality
\[
 \genfrac{}{}{}{}{\dim\Hom((L_m)_k,F)}{\dim\Hom((L_n)_k,F)}<
 \genfrac{}{}{}{}{\dim\Hom((L_m)_k,E)}{\dim\Hom((L_n)_k,E)}
 \quad
 \left(
 \text{resp.} \;
 \genfrac{}{}{}{}{\dim\Hom((L_m)_k,F)}{\dim\Hom((L_n)_k,F)} \leq
 \genfrac{}{}{}{}{\dim\Hom((L_m)_k,E)}{\dim\Hom((L_n)_k,E)}
 \right)
\]
holds for $n\gg m\gg 0$ and for any non-zero object
$F\in{\mathcal D}_k$ satisfying
$\Ext^i((L_N)_k,F)=0$ for $N\gg 0$ and $i\neq 0$ with a morphism
$\iota:F\rightarrow E$ such that $\iota$ is not isomorphic and
$\Hom((L_n)_k,F)\rightarrow\Hom((L_n)_k,E)$ is injective for $n\gg 0$.
\end{definition}

\begin{remark}\label{stability-condition}\rm
Let $\Spec k\rightarrow S$ be a geometric point
and $E$ an object of ${\mathcal D}_k$ satisfying
$\Ext^i((L_n)_k,E)=0$ for $i\neq 0$ and $n\gg 0$.
Let ${\mathcal E}$ be the coherent ${\mathcal A}_k$-module corresponding
to $E$ as in Proposition \ref{equivalence}.
Then $E$ is ${\mathcal L}$-stable (resp.\ ${\mathcal L}$-semistable)
if and only if for any coherent ${\mathcal A}_k$-submodule ${\mathcal F}$
of ${\mathcal E}$ with $0\neq{\mathcal F}\subsetneq{\mathcal E}$,
the inequality
\begin{equation}\label{A-mod-stability}
 \genfrac{}{}{}{}{\chi({\mathcal F}(m))}{\chi({\mathcal F}(n))} <
 \genfrac{}{}{}{}{\chi({\mathcal E}(m))}{\chi({\mathcal E}(n))}
 \quad
 \left( \text{resp.} \;
 \genfrac{}{}{}{}{\chi({\mathcal F}(m))}{\chi({\mathcal F}(n))} \leq
 \genfrac{}{}{}{}{\chi({\mathcal E}(m))}{\chi({\mathcal E}(n))}
 \right)
\end{equation}
holds for $n\gg m\gg0$.
We say a coherent ${\mathcal A}_k$-module ${\mathcal E}$
stable (resp.\ semistable) if the corresponding object
$E$ of ${\mathcal D}_k$ is ${\mathcal L}$-stable
(resp.\ ${\mathcal L}$-semistable).
\end{remark}

\begin{remark}\rm
For a field $K$ with a morphism $\Spec K \rightarrow S$
and an object $E\in{\mathcal D}_K$,
we say that $E$ is ${\mathcal L}$-stable (resp.\ ${\mathcal L}$-semistable) if
$E_{\bar{K}}$ is ${\mathcal L}$-stable (resp.\ ${\mathcal L}$-semistable),
where $\bar{K}$ is the algebraic closure of $K$.
\end{remark}

\section{Existence of the moduli space of stable objects}

\begin{definition}\rm
 Let $p:{\mathcal D}\rightarrow(\sch/S)$ be a fibered triangulated category
 with base change property and ${\mathcal L}=\{L_n\}_{n\geq 0}$
 be a strict ample sequence.
 For a numerical polynomial $P(t)\in\mathbf{Q}[t]$, we define a moduli functor
 ${\mathcal M}_{\mathcal D}^{P,\mathcal L}:(\sch/S)\rightarrow(\sets)$ by
\[
 {\mathcal M}_{\mathcal D}^{P,\mathcal L}(T):=
 \left\{
 E\in{\mathcal D}_T \left|
 \begin{array}{l}
  \text{for any geometric point $s$ of $T$, for $n\gg 0$,
  $\Ext^i((L_n)_s,E_s)=0$} \\
  \text{for $i\neq 0$ and $\Hom((L_n)_s,E_s)=P(n)$ and
  $E_s$ is ${\mathcal L}$-stable}
 \end{array}
 \right\}\right/\sim,
\]
where $E\sim E'$ if there exists a line bundle $L$ on $T$
and an isomorphism $E\stackrel{\sim}\rightarrow E'\otimes L$.

We also define a moduli functor
$\overline{\mathcal M}_{\mathcal D}^{P,\mathcal L}:
(\sch/S)\rightarrow(\sets)$ by
\[
 \overline{{\mathcal M}_{\mathcal D}^{P,\mathcal L}}(T):=
 \left\{
 E\in{\mathcal D}_T \left|
 \begin{array}{l}
  \text{for any geometric point $s$ of $T$, for $n\gg 0$,
  $\Ext^i((L_n)_s,E_s)=0$ } \\
  \text{for $i\neq 0$ and $\Hom((L_n)_s,E_s)=P(n)$ and
  $E_s$ is ${\mathcal L}$-semistable}
 \end{array}
 \right\}\right/\sim,
\]
where $E\sim E'$ if there exists a line bundle $L$ on $T$
such that $E\cong E'\otimes L$ or
there exist sequences
$0=E_0\rightarrow E_1\rightarrow\cdots\rightarrow E_{\alpha}=E$
and
$0=E'_0\rightarrow E'_1\rightarrow\cdots\rightarrow E'_{\alpha}=E'$
such that 
$\Ext^i((L_n)_s,(E_j)_s)=\Ext^i((L_n)_s,(E'_j)_s)=0$
for $n\gg 0$, $i\neq 0$ and $s\in T$,
$\Hom((L_n)_s,(E_j)_s)\rightarrow \Hom((L_n)_s,(E_{j+1})_s)$ and
$\Hom((L_n)_s,(E'_j)_s)\rightarrow \Hom((L_n)_s,(E'_{j+1})_s)$
are injective for $n\gg 0$ and $s\in T$ and
$\bigoplus_{j=1}^{\alpha} F_j\cong\bigoplus_{j=1}^{\alpha} F'_j\otimes L$,
where $F_j=\mathrm{Cone}(E_{j-1}\rightarrow E_j)$,
$F'_j=\mathrm{Cone}(E'_{j-1}\rightarrow E'_j)$
and for any geometric point $s$ of $T$,
$(F_j)_s$ and $(F'_j)_s$ are ${\mathcal L}$-stable
such that
\[
 \genfrac{}{}{}{}{\dim\Hom((L_m)_s,(F_j)_s)}{\dim\Hom((L_n)_s,(F_j)_s)}
 =\genfrac{}{}{}{}{P(m)}{P(n)}
 =\genfrac{}{}{}{}{\dim\Hom((L_m)_s,(F'_j)_s)}{\dim\Hom((L_n)_s,(F'_j)_s)}
\]
for $n\gg m\gg 0$ and for $j=1,2,\ldots,\alpha$.
\end{definition}

\begin{proposition}\label{boundedness}
 For any numerical polynomial $P(t)\in\mathbf{Q}[t]$, the family
\[
 \left\{ E \left|
 \begin{array}{l}
 \text{\rm $E\in{\mathcal D}_k$ for some geometric point
 $\Spec k\rightarrow S$,}\\
 \text{\rm $E$ is ${\mathcal L}$-semistable and
 $\Hom((L_n)_k,E)=P(n)$ for $n\gg 0$}
 \end{array}
 \right\}\right.
\]
is bounded.
\end{proposition}

\begin{proof}
It suffices to show that the corresponding family of
coherent ${\mathcal A}$-modules on the fibers of $\mathbf{P}(V_1)$
over $S$ is bounded.
For a coherent sheaf ${\mathcal G}$ on $\mathbf{P}(V_1)$, we can write
\[
 \chi({\mathcal G}(n))=\sum_{i=0}^d a_i({\mathcal G})\genfrac{(}{)}{0pt}{}{n+d-i}{d-i}
\]
with $a_i({\mathcal G})$ integers
and we write $\mu(G)=a_1({\mathcal G})/a_0({\mathcal G})$.
Let ${\mathcal E}$ be a coherent ${\mathcal A}_k$-module such that
$\chi({\mathcal E}(n))=P(n)$ and
the corresponding object of ${\mathcal D}_k$ is ${\mathcal L}$-semistable.
Note that ${\mathcal E}$ is of pure dimension.
We can take the slope maximal destabilizer ${\mathcal F}$ of ${\mathcal E}$
as a sheaf on $\mathbf{P}(V_1)$.
Let $\tilde{\mathcal F}$ be the image of
${\mathcal F}\otimes{\mathcal A}\rightarrow {\mathcal E}$.
Note that there exists a locally free sheaf $W$ of finite rank on $S$, positive integer $N$
and a surjection
\[
 W\otimes {\mathcal O}(-N)\longrightarrow{\mathcal A}
\]
Then we obtain a surjection
\[
 W\otimes{\mathcal F}(-N)\longrightarrow {\mathcal F}\otimes{\mathcal A}\longrightarrow \tilde{\mathcal F}.
\]
Since $W\otimes{\mathcal F}(-N)$ is slope semistable, we have
\[
 \mu({\mathcal F})-N=\mu(W\otimes{\mathcal F}(-N))
 \leq \mu(\tilde{\mathcal F})\leq\mu({\mathcal E}).
\]
So the maximal slope $\mu({\mathcal F})$ is bounded by $N+\mu({\mathcal E})$.
Then we obtain the boundedness by [\cite{langer}, Theorem 4.2].
\end{proof}

\begin{proposition}\label{openness-stability}
 Assume that $U\in(\sch/S)$ and $E\in{\mathcal D}_U$ are given.
 Then the subsets
 \begin{align*}
  U^s &= \left\{ x\in U \left| 
  \text{$E_x$ is ${\mathcal L}$-stable}\right\}\right. \\
  U^{ss} &= \left\{ x\in U \left|
  \text{$E_x$ is ${\mathcal L}$-semistable}\right\}\right.
 \end{align*}
 of $U$ are open.
\end{proposition}

\begin{proof}
First we will show that
\[
 U'=\left\{ x\in U \left|
 \text{$\Ext^i((L_n)_x,E_x)=0$ for $n\gg 0$ and $i\neq 0$}
 \right\} \right.
\]
is open in $U$.
By Definition \ref{def-polarization} (4),
there exists a positive integer $n_0$ such that for any $n\geq n_0$,
exists an integer $N_n$ with $N_n>n$
such that for any $N\geq N_n$,
\[
 \Hom((L_N)_s,(L_n)_s)\otimes \Ext^i((L_n)_s,E_s)
 \longrightarrow \Ext^i((L_N)_s,E_s)
\]
is surjective for any $i$ and $s\in U$.
By Definition \ref{def-polarization} (4),
there are integers $k_1,k_2$ with $k_1<k_2$ such that
$\Ext^i((L_{n_0})_s,E_s)=0$ for any $s\in U$ except for $k_1\leq i\leq k_2$.
Then we have
$\Ext^i((L_N)_s,E_s)=0$ for $N\geq N_{n_0}$ and $s\in U$,
except for $k_1\leq i\leq k_2$.
Now take any point $x\in U'$.
For each $i\neq 0$ with $k_1\leq i\leq k_2$, there is an integer
$m_i$ with $m_i\geq n_0$ such that $\Ext^i((L_{m_i})_x,E_x)=0$.
For any $N\geq N_{m_i}$,
\[
 \Hom((L_N)_s,(L_{m_i})_s)\otimes\Ext^i((L_{m_i})_s,E_s)
 \longrightarrow \Hom((L_N)_s,E_s)
\]
is surjective for any $s\in U$.
By using Definition \ref{bcp} (4), we can see that there exists an open neighborhood $U_i$ of $x$ such that
$\Ext^i((L_{m_i})_y,E_y)=0$ for any $y\in U_i$.
Then we have $\Ext^i((L_N)_y,E_y)=0$ for $N\geq N_{m_i}$.
If we put
\[
 V:=\bigcap_{k_1\leq i\leq k_2, i\neq 0} U_i
\]
then $V$ is an open neighborhood of $x$.
Put
\[
 \tilde{N}:=\max (\{ N_{m_i}|k_1\leq i\leq k_2, i\neq 0\}\cup\{N_{n_0}\}).
\]
Then we have $\Ext^i((L_N)_y,E_y)=0$
for any $y\in V$, $i\neq 0$ and $N\geq \tilde{N}$,
which means $V\subset U'$.
Thus $U'$ is an open subset of $U$.

By Proposition \ref{equivalence}, $E_{U'}$ corresponds to a coherent ${\mathcal A}_{U'}$-module
${\mathcal E}$ flat over $U'$.
We can see that $U^s$ coincides with
\[
 \left\{ x\in U' \left| \text{${\mathcal E}\otimes k(x)$ is a stable
 ${\mathcal A}_x$-module} \right\}\right..
\]
We can see by the argument similar to that of [\cite{H-L}, Proposition 2.3.1]
that this subset is open in $U'$.
By the same argument we can also see the openness of $U^{ss}$.
\end{proof}

\begin{theorem}\label{thm:moduli-exists}
 There exists a coarse moduli scheme
 $\overline{M_{\mathcal D}^{P,\mathcal L}}$
 of $\overline{{\mathcal M}_{\mathcal D}^{P,\mathcal L}}$
 and an open subscheme $M_{\mathcal D}^{P,\mathcal L}$
 of $\overline{M_{\mathcal D}^{P,\mathcal L}}$
 which is a coarse moduli scheme of
 ${\mathcal M}_{\mathcal D}^{P,\mathcal L}$.
\end{theorem}

Before constructing the moduli space, we first note the following lemma:

\begin{lemma}\label{fundamental-lemma}
 Let $P(x)$ be a numerical polynomial.
 Then there exists an integer $m_0$ such that for any $m\geq m_0$,
 any geometric point $s$ of $S$, any semi-stable ${\mathcal A}_s$-module
 ${\mathcal E}$ with $\chi({\mathcal E}(n))=P(n)$,
\begin{enumerate}
\item ${\mathcal E}(m)$ is generated by global sections and
 $H^i({\mathcal E}(m))=0$ for $i>0$,
\item for any nonzero coherent ${\mathcal A}_s$-submodule
 ${\mathcal F}\subset{\mathcal E}$, the inequality
\[
 \dim H^0({\mathcal F}(m))\leq
 \frac{a_0({\mathcal F})}{a_0({\mathcal E})}
 \dim H^0({\mathcal E}(m))
\]
holds, where
\[
 \chi({\mathcal E}(n))=\sum_{i=0}^d a_i({\mathcal E})
 \genfrac{(}{)}{0pt}{}{n+d-i}{d-i},
 \quad
 \chi({\mathcal F}(n))=\sum_{i=0}^d a_i({\mathcal F})
 \genfrac{(}{)}{0pt}{}{n+d-i}{d-i}.
\]
Moreover the equality holds if and only if
$\chi({\mathcal E}(n))/a_0({\mathcal E})=
\chi({\mathcal F}(n))/a_0({\mathcal F})$
as polynomials in $n$.
\end{enumerate} 
\end{lemma}

\begin{proof}
Proof is essentially the same as [\cite{maruyama}, Proposition 4.10.]
\end{proof}

Take $m_0$ as in Lemma \ref{fundamental-lemma}.
Replacing $S$ by its connected component,
we may assume that $S$ is connected.
Replacing $m_0$ if necessary, we may assume by Proposition \ref{boundedness}
that for any geometric point
$E\in \overline{{\mathcal M}_{\mathcal D}^{P,\mathcal L}}(k)$
and for any $m\geq m_0$,
$\Ext^i((L_m)_k,E)=0$ for $i\neq 0$ and
\[
 \Hom((L_n)_k,(L_m)_k)\otimes\Hom((L_m)_k,E)
 \longrightarrow\Hom((L_n)_k,E)
\]
is surjective for $n\gg 0$.
For a geometric point
$E\in\overline{{\mathcal M}_{\mathcal D}^{P,\mathcal L}}(k)$,
we consider the canonical morphism
\[
 u:(L_{m_0})_k\otimes\Hom((L_{m_0})_k,E)\longrightarrow E
\]
and put $E_1:=\mathrm{Cone}(u)[-1]$.
We can take $m_1\gg m_0$ such that for any such $E$ and
for any $m\geq m_1$, $\Ext^i((L_m)_k,E_1)=0$ for $i\neq 0$ and
\[
 \Hom((L_n)_k,(L_m)_k)\otimes\Hom((L_m)_k,E_1)
 \longrightarrow \Hom((L_n)_k,E_1)
\]
is surjective for $n\gg 0$.
We consider the canonical morphism
\[
 v:(L_{m_1})_k\otimes\Hom((L_{m_1})_k,E_1)\longrightarrow E_1
\]
and put $E_2:=\mathrm{Cone}(v)[-1]$.
We can take $m_2\gg 0$ such that for any $E$ and for any $m\geq m_2$,
$\Ext^i((L_m)_k,E_2)=0$ for $i\neq 0$ and
\[
 \Hom((L_n)_k,(L_m)_k)\otimes\Hom((L_m)_k,E_2)
 \longrightarrow\Hom((L_n)_k,E_2)
\]
is surjective for $n\gg 0$.
We put
\[
 r_0:=\dim_k\Hom((L_{m_0})_k,E), \quad
 r_1:=\dim_k((L_{m_1})_k,E_1), \quad
 r_2:=\dim_k((L_{m_2})_k,E_2)
\]
and
\[
 W_0:={\mathcal O}_S^{\oplus r_0}, \quad
 W_1:={\mathcal O}_S^{\oplus r_1}, \quad
 W_2:={\mathcal O}_S^{\oplus r_2}.
\]
Note that $r_0,r_1,r_2$ are independent of the choice of $E$
and only depend on $P$ and ${\mathcal L}$.
We set
\[
 Z:=
 \mathbf{V}\left(R^0\Hom_p(L_{m_2},L_{m_1})^{\vee}
 \otimes W_2\otimes W_1^{\vee}\right)
 \times
 \mathbf{V}\left(R^0\Hom_p(L_{m_1},L_{m_0})^{\vee}
 \otimes W_1\otimes W_0^{\vee}\right).
\]
Let
\[
 (L_{m_2})_Z\otimes W_2 \stackrel{\tilde{v}}\longrightarrow
 (L_{m_1})_Z\otimes W_1 \stackrel{\tilde{u}}\longrightarrow
 (L_{m_0})_Z\otimes W_0
\]
be the universal family.
There exists a closed subscheme $Y\subset Z$ such that
\[
 Y(T)=\left\{ g\in Z(T) \left|
 g^*(\tilde{u}\circ \tilde{v})=0 \right\}\right.
\]
for any $T\in(\sch/S)$.
Since the sequence
\begin{gather*}
 \Hom(\mathrm{Cone}(\tilde{v}_Y),(L_{m_0})_Y\otimes W_0)
 \stackrel{\beta}\longrightarrow
 \Hom((L_{m_1})_Y\otimes W_1,(L_{m_0})_Y\otimes W_0) \\
 \stackrel{\tilde{v}^*}\longrightarrow
 \Hom((L_{m_2})_Y\otimes W_2,(L_{m_0})_Y\otimes W_0)
\end{gather*}
is exact and $\tilde{v}^*(\tilde{u}_Y)=\tilde{u}_Y\circ\tilde{v}_Y=0$,
there exists a morphism
$\tilde{w}:\mathrm{Cone}(\tilde{v}_Y)\rightarrow (L_{m_0})_Y\otimes W_0$
such that $\beta(\tilde{w})=\tilde{u}_Y$.
We put $\tilde{B}:=\mathrm{Cone}(\tilde{w})$ and set
\[
 Y':=\left\{ x\in Y \left|
  \text{$\Ext^{-1}((L_n)_x,\tilde{B}_x)=0$ for $n\gg 0$}
 \right\}\right.
\]
Then we can see that $Y'$ is an open subset of $Y$.
Note that for any $x\in Y'$,
$\Ext^i((L_n)_x,\tilde{B}_x)=0$ for $n\gg 0$
except for $i=-2,0$.
By Definition \ref{def-polarization} (5),
there exist an object $\tilde{E}\in{\mathcal D}_{Y'}$
and a morphism $\tilde{B}_{Y'}\rightarrow\tilde{E}$ such that
$\Ext^i((L_n)_x,\tilde{E}_x)=0$ for $n\gg 0$, $x\in Y'$ and $i\neq 0$
and $\Hom((L_n)_x,\tilde{B}_x)\rightarrow\Hom((L_n)_x,\tilde{E}_x)$
is isomorphic for $n\gg 0$ and $x\in Y'$.
If we set
\[
 \tilde{E}_1:=\mathrm{Cone}
 ((L_{m_0})_{Y'}\otimes W_0\rightarrow\tilde{E})[-1],
\]
$\mathrm{Cone}(\tilde{v})_{Y'}\rightarrow (L_{m_0})_{Y'}\otimes W_0$
factors through $\tilde{E}_1$.
Moreover, for any $x\in Y'$,
$\Ext^i((L_n)_x,(\tilde{E}_1)_x)=0$ for $i\neq 0$ and
$\Hom((L_n)_x,\mathrm{Cone}(\tilde{v})_x)
\rightarrow\Hom((L_n)_x,(\tilde{E}_1)_x)$
is isomorphic for $n\gg 0$.
If we set
\[
 \tilde{E}_2:=\mathrm{Cone}
 ((L_{m_1})_{Y'}\otimes W_1\rightarrow\tilde{E}_1)[-1],
\]
then $\tilde{v}_{Y'}$ factors through $\tilde{E}_2$.
Now we put
\[
 Y^{ss}:=\left\{ x\in Y' \left|
 \begin{array}{l}
  \text{$W_0\otimes k(x)\rightarrow \Hom((L_{m_0})_x,\tilde{E}_x)$
  is isomorphic,} \\
  \text{$W_j\otimes k(x)\rightarrow \Hom((L_{m_j})_x,(\tilde{E}_j)_x)$
  are isomorphic for $j=1,2$,} \\
  \text{$\Hom((L_n)_x,\tilde{E}_x)=P(n)$ for $n\gg 0$
  and $\tilde{E}_x$ is ${\mathcal L}$-semistable}
 \end{array}
 \right\}\right.
\]
and
\[
 Y^s:=\left\{ x\in Y^{ss} \left|
 \text{$\tilde{E}_x$ is ${\mathcal L}$-stable}
 \right\}\right..
\]
Then we can check that $Y^s, Y^{ss}$ are open subsets of $Y'$.
If we put
\[
 G:=GL(W_0)\times GL(W_1)\times GL(W_2),
\]
then there is a canonical action of $G$ on $Z$
and $Y$, $Y'$, $Y^{ss}$, $Y^s$ are preserved by this action.
For a sufficiently large integer $N$, we put
\begin{align*}
 \alpha_0&:=\rk W_2+N\rk W_1 \\
 \alpha_1&:=-N\rk W_0 \\
 \alpha_2&:=-\rk W_0
\end{align*}
and consider the character
\[
 \chi: G \longrightarrow \mathbf{G}_m ; \quad
 (g_0,g_1,g_2) \mapsto
 \det(g_0)^{\alpha_0}\det(g_1)^{\alpha_1}\det(g_2)^{\alpha_2}.
\]
Let us consider the quiver consisting of three vertices
$v_2,v_1,v_0$ and
$\rk_{{\mathcal O}_S}R^0\Hom_p(L_{m_2},L_{m_1})$-arrows
from $v_2$ to $v_1$ and
$\rk_{{\mathcal O}_S}R^0\Hom_p(L_{m_1},L_{m_0})$-arrows
from $v_1$ to $v_0$.
Then the points of $Z$ correspond to the representations
of this quiver (see \cite{king} for the definition of
quiver and its representation).

\begin{lemma}\label{lemm:GIT-stable}
 If we take $N\gg m_2\gg m_1\gg m_0\gg 0$,
 $Y^{ss}$ is contained in the set
 $Z^{ss}(\chi)$ of $\chi$-semistable points of $Z$ 
 in the sense of \cite{king}.
 Moreover, $Y^s$ is contained in the set $Z^s(\chi)$
 of $\chi$-stable points of $Z$.
\end{lemma}

\begin{proof}
Take any geometric point $x$ of $Y^{ss}$
and vector subspaces
$W'_i\subset (W_i)_x$ ($0\leq i\leq 2$)
which induce commutative diagrams
\[
\begin{CD}
 W'_2 @>>> W'_1\otimes R^0\Hom_p(L_{m_2},L_{m_1})_x \\
 @VVV @VVV \\
(W_2)_x @>>> (W_1)_x\otimes R^0\Hom_p(L_{m_2},L_{m_1})_x
\end{CD}
\quad
\begin{CD}
 W'_1 @>>> W'_0\otimes R^0\Hom_p(L_{m_1},L_{m_0})_x \\
 @VVV @VVV \\
 (W_1)_x @>>> (W_0)_x\otimes R^0\Hom_p(L_{m_1},L_{m_0})_x.
\end{CD}
\]
From \cite{king}, we should say that
\[
 \alpha_0\dim W'_0 + \alpha_1\dim W'_1+\alpha_2\dim W'_2\geq 0.
\]
Let ${\mathcal E}$ be the $Y^{ss}$-flat ${\mathcal A}_{Y^{ss}}$-module
corresponding to $\tilde{E}_{Y^{ss}}$ by Proposition \ref{equivalence}.
Then a morphism
${\mathcal A}(-m_0)\otimes W'_0\rightarrow {\mathcal E}_x$
is induced and we denote its image by ${\mathcal E}(W'_0)$.
Note that ${\mathcal E}_x$ is of pure dimension and so
${\mathcal E}(W'_0)$ is also of pure dimension.
Since the family
\[
 \left\{ {\mathcal E}(W'_0) \left|
 W'_0\subset (W_0)_x,\; \text{$x$ is a geometric point of $Y^{ss}$} \right\}\right.
\]
is bounded, we can find an integer $m_1\gg m_0$ such that for
$K'_1:=\ker(W'_0\otimes{\mathcal A}(-m_0)\rightarrow{\mathcal E}(W'_0))$,
$K'_1(m_1)$ is generated by global sections and
$H^i(K'_1(m_1))=0$,  $H^i({\mathcal A}_x(m_1-m_0))=0$ for $i>0$.
Moreover we can find an integer $m_2\gg m_1$ such that for
$K'_2:=\ker(H^0(K'_1(m_1))\otimes{\mathcal A}(-m_1)\rightarrow K'_1)$,
$K'_2(m_2)$ is generated by global sections and
$H^i(K'_2(m_2))=0$, $H^i({\mathcal A}_x(m_2-m_1))=0$, $H^i({\mathcal A}_x(m_2-m_0))=0$
and $H^i(K'_1(m_2))=0$ for $i>0$.
If we put $\tilde{W}'_1:=H^0(K'_1(m_1))$ and
$\tilde{W}'_2:=H^0(K'_2(m_2))$, then we have
\begin{align*}
 \dim H^0({\mathcal E}(W'_0)(m_1)) &=
 \dim H^0({\mathcal A}_x(m_1-m_0))\dim W'_0-\dim \tilde{W}'_1 \\
 \dim H^0({\mathcal E}(W'_0)(m_2)) &=
 \dim H^0({\mathcal A}_x(m_2-m_0))\dim W'_0
 -\dim H^0({\mathcal A}_x(m_2-m_1))\dim \tilde{W}'_1 +\dim \tilde{W}'_2.
\end{align*}
Since the family $\{{\mathcal  E}(W'_0)\}$ is bounded, we can take by using Lemma \ref{fundamental-lemma}
a positive integer $m_0\gg 0$
and a positive number $\epsilon>0$ such that
\[
 \genfrac{}{}{}{}{h^0({\mathcal E}(W'_0)(m_0))}{P(m_0)}<\genfrac{}{}{}{}{a_0({\mathcal E}(W'_0))}{a_0(P)}-\epsilon
\]
for any $W'_0$ such that
\[
 \genfrac{}{}{}{}{\chi({\mathcal E}(W'_0)(m))}{\chi({\mathcal E}(W'_0)(n))}<
 \genfrac{}{}{}{}{P(m)}{P(n)}
\]
for $n\gg m\gg 0$.
Here we write 
\[
 \chi({\mathcal E}(W'_0)(n))=\sum_{i=0}^d a_i({\mathcal E}(W'_0)) \genfrac{(}{)}{0pt}{}{n+d-i}{d-i},
 \quad
 P(n)=\sum_{i=0}^d a_i(P) \genfrac{(}{)}{0pt}{}{n+d-i}{d-i}
\]
with $a_i({\mathcal E}(W'_0))$ and $a_i(P)$ integers.
Since 
\[
 \lim_{m_1\to\infty}\genfrac{}{}{}{}{h^0({\mathcal E}(W'_0)(m_1))}{P(m_1)}=
 \genfrac{}{}{}{}{a_0({\mathcal E}(W'_0))}{a_0(P)},
\]
we can take $m_1\gg m_0$ such that
\[
 \genfrac{}{}{}{}{h^0({\mathcal E}(W'_0)(m_1))}{P(m_1)} >
 \genfrac{}{}{}{}{a_0({\mathcal E}(W'_0))}{a_0(P)}-\genfrac{}{}{}{}{\epsilon}{2}.
\]
Since
\[
 \lim_{N\to\infty}\genfrac{}{}{}{}
 {(h^0({\mathcal A}_x(m_2-m_1))+N)h^0({\mathcal E}(W'_0)(m_1))
 -h^0({\mathcal E}(W'_0)(m_2))}
 {(h^0({\mathcal A}_x(m_2-m_1))+N)P(m_1)-P(m_2)}=
 \genfrac{}{}{}{}{h^0({\mathcal E}(W'_0)(m_1))}{P(m_1)},
 \]
we can take $N\gg m_2$ such that
 \[
  \genfrac{}{}{}{}{(h^0({\mathcal A}_x(m_2-m_1))+N)h^0({\mathcal E}(W'_0)(m_1))
  -h^0({\mathcal E}(W'_0)(m_2))}
  {(h^0({\mathcal A}_x(m_2-m_1))+N)P(m_1)-P(m_2)}
  >
  \genfrac{}{}{}{}{h^0({\mathcal E}(W'_0)(m_1))}{P(m_1)}-\genfrac{}{}{}{}{\epsilon}{2}.
 \]
Then we have
 \begin{align*}
  \genfrac{}{}{}{}{h^0({\mathcal E}(W'_0)(m_0))}{P(m_0)} &<
  \genfrac{}{}{}{}{a_0({\mathcal E}(W'_0))}{a_0(P)}-\epsilon \\
  &< \genfrac{}{}{}{}{h^0({\mathcal E}(W'_0)(m_1))}{P(m_1)}+\genfrac{}{}{}{}{\epsilon}{2}-\epsilon \\
  &< \genfrac{}{}{}{}{(h^0({\mathcal A}_x(m_2-m_1))+N)h^0({\mathcal E}(W'_0)(m_1))
  -h^0({\mathcal E}(W'_0)(m_2))}
  {(h^0({\mathcal A}_x(m_2-m_1))+N)P(m_1)-P(m_2)}+\genfrac{}{}{}{}{\epsilon}{2}+\genfrac{}{}{}{}{\epsilon}{2}-\epsilon \\
  &= \genfrac{}{}{}{}{(h^0({\mathcal A}_x(m_2-m_1))+N)h^0({\mathcal E}(W'_0)(m_1))
  -h^0({\mathcal E}(W'_0)(m_2))}
  {(h^0({\mathcal A}_x(m_2-m_1))+N)P(m_1)-P(m_2)}
 \end{align*}
for any $W'_0$ such that
 \[
  \genfrac{}{}{}{}{\chi({\mathcal E}(W'_0)(m))}{\chi({\mathcal E}(W'_0)(n))}
  <\genfrac{}{}{}{}{P(m)}{P(n)}
 \]
for $n\gg m\gg 0$.
Take $W'_0$ such that
 \[
  \genfrac{}{}{}{}{\chi({\mathcal E}(W'_0)(m))}{\chi({\mathcal E}(W'_0)(n)}=
  \genfrac{}{}{}{}{P(m)}{P(n)} 
 \]
for $n\gg m\gg 0$.
Then we can see by Lemma \ref{fundamental-lemma} that
 \[
  \genfrac{}{}{}{}{h^0({\mathcal E}(W'_0)(m_0))}{P(m_0)}=\genfrac{}{}{}{}{a_0({\mathcal E}(W'_0))}{a_0(P)}
  = \genfrac{}{}{}{}{(h^0({\mathcal A}_x(m_2-m_1))+N)h^0({\mathcal E}(W'_0)(m_1))
  -h^0({\mathcal E}(W'_0)(m_2))}
   {(h^0({\mathcal A}_x(m_2-m_1))+N)P(m_1)-P(m_2)}.
 \]
Hence we have the inequality
\begin{equation}\label{GIT-stability-inequality}
 h^0({\mathcal E}(W'_0)(m_0))\leq
 \genfrac{}{}{}{}
 {(h^0({\mathcal A}_x(m_2-m_1))+N)h^0({\mathcal E}(W'_0)(m_1))
 -h^0({\mathcal E}(W'_0)(m_2)}
 {(h^0({\mathcal A}_x(m_2-m_1))+N)P(m_1)-P(m_2)}P(m_0)
\end{equation}
for any ${\mathcal E}(W'_0)$.
Moreover, the equality holds in (\ref{GIT-stability-inequality}) if and only if
$\chi({\mathcal E}(W'_0)(n))/a_0({\mathcal E}(W'_0))=P(n)/a_0(P)$
as polynomials in $n$.
From the inequality (\ref{GIT-stability-inequality}), we obtain the inequality
\[
 (r_2+Nr_1)\dim W'_0-Nr_0\dim \tilde{W}'_1-r_0\dim \tilde{W}'_2\geq 0
\]
by using $\dim W'_0\leq h^0({\mathcal E}(W'_0)(m_0))$.
Since $\dim W'_1\leq \dim \tilde{W}'_1$ and $\dim W'_2\leq\dim\tilde{W}'_2$,
we have
\begin{equation}\label{quiver-inequality}
 \alpha_0\dim W'_0+\alpha_1\dim W'_1+\alpha_2\dim W'_2\geq 0.
\end{equation}
Thus $x$ becomes a geometric point of $Z^{ss}(\chi)$.

In the inequality (\ref{quiver-inequality}), the equality holds
if and only if $\dim \tilde{W}'_1=\dim W'_1$, $\dim \tilde{W}'_2=\dim W'_2$,
$h^0({\mathcal E}(W'_0))=\dim W'_0$ and
$\chi({\mathcal E}(W'_0)(n))/a_0({\mathcal E}(W'_0))=P(n)/a_0(P)$
as polynomials in $n$.
So, if $x$ is a geometric point of $Y^s$, we have
\[
 (r_2+Nr_1)\dim W'_0-Nr_0\dim W'_1-r_0\dim W'_2 > 0.
\]
for any $(W'_0,W'_1,W'_2)$ with $(0,0,0)\neq(W'_0,W'_1,W'_2)\subsetneq((W_0)_x,(W_1)_x,(W_2)_x)$,
which means that $x$ becomes a geometric point of $Z^s(\chi)$.
\end{proof}

By \cite{king} and \cite{seshadri}, there exists a GIT quotient
$\phi:Y\cap Z^{ss}(\chi)\rightarrow (Y\cap Z^{ss}(\chi))//G$.

\begin{lemma}\label{open-lemma}
 $\phi^{-1}(\phi(Y^{ss}))=Y^{ss}$.
\end{lemma}

\begin{proof}
It is sufficient to show that $\phi^{-1}(\phi(Y^{ss}))\subset Y^{ss}$.
Take any $k$-valued geometric point $x$ of $\phi^{-1}(\phi(Y^{ss}))$.
Let $s$ be the induced $k$-valued goemetric point of $S$.
Since $\phi(x)$ is a geometric point of $\phi(Y^{ss})$,
there exists a $k$-valued geometric point $y$ of $Y^{ss}$
such that $\phi(x)=\phi(y)$.
         
Let ${\mathcal E}$ be the $Y^{ss}$-flat ${\mathcal A}_{Y^{ss}}$-module corresponding to $\tilde{E}_{Y^{ss}}$
as in the proof of Lemma \ref{lemm:GIT-stable}.
Then there is a Jordan-H\"older filtration
\[
 0=F^{(0)}\subset F^{(1)}\subset\cdots\subset F^{(l)}={\mathcal E}\otimes k(y)
\]
of ${\mathcal E}\otimes k(y)$.
For each $i$ with $1\leq i\leq l$, we define $K^{(i)}_1$, $K^{(i)}_2$
by exact sequences
\begin{gather*}
 0\longrightarrow K^{(i)}_1 \longrightarrow
 H^0(F^{(i)}(m_0))\otimes{\mathcal A}(-m_0)
 \longrightarrow F^{(i)} \longrightarrow 0 \\
 0\longrightarrow K^{(i)}_2 \longrightarrow
 H^0(K^{(i)}_1(m_1))\otimes{\mathcal A}(-m_1)
 \longrightarrow K^{(i)}_1 \longrightarrow 0.
\end{gather*}
Then $y$ corresponds to the representation of quiver given by
\begin{gather*}
 H^0(K^{(l)}_2(m_2))\longrightarrow
 H^0(K^{(l)}_1(m_1))\otimes H^0({\mathcal A}_s(m_2-m_1)) \\
 H^0(K^{(l)}_1(m_1))\longrightarrow
 H^0(F^{(l)}(m_0))\otimes H^0({\mathcal A}_s(m_1-m_0))
\end{gather*}
and the Jordan-H\"older filtration of ${\mathcal E}\otimes k(y)$
corresponds to the filtration of the quiver representation given by
\begin{gather*}
 0\subset H^0(K^{(1)}_2(m_2))\subset\cdots\subset H^0(K^{(l)}_2(m_2)) \\
 0\subset H^0(K^{(1)}_1(m_1))\subset\cdots\subset H^0(K^{(l)}_1(m_1)) \\
 0\subset H^0(F^{(1)}(m_0))\subset\cdots\subset H^0(F^{(l)}(m_0)).
\end{gather*}

We put $E^{(i)}:=F^{(i)}/F^{(i-1)}$ and
$\overline{\mathcal E}:=\bigoplus_{i=1}^l E^{(i)}$.
For $i=1,\ldots,l$, we define $\bar{K}^{(i)}_1$, $\bar{K}^{(i)}_2$
by the exact sequences
\begin{gather*}
 0\longrightarrow \bar{K}^{(i)}_1 \longrightarrow
 H^0(E^{(i)}(m_0))\otimes{\mathcal A}(-m_0) \longrightarrow
 E^{(i)} \longrightarrow 0 \\
 0\longrightarrow \bar{K}^{(i)}_2 \longrightarrow
 H^0(\bar{K}_1^{(i)}(m_1))\otimes{\mathcal A}(-m_1) \longrightarrow
 \bar{K}^{(i)}_1 \longrightarrow 0. \\
\end{gather*}
We can see from the proof of Lemma \ref{lemm:GIT-stable}
that the quiver representation $y_i$ given by
\begin{gather*}
 H^0(\bar{K}^{(i)}_2(m_2))\longrightarrow
 H^0(\bar{K}^{(i)}_1(m_1))\otimes H^0({\mathcal A}_s(m_2-m_1)) \\
 H^0(\bar{K}^{(i)}_1(m_1))\longrightarrow
 H^0(E^{(i)}(m_0))\otimes H^0({\mathcal A}_s(m_1-m_0))
\end{gather*}
is stable with respect to the weight $(\alpha_0,\alpha_1,\alpha_2)$.
The direct sum $y_1\oplus\cdots\oplus y_l$ corresponds
to a point $y'$ of $Y^{ss}_s$ given by the exact sequence
\begin{gather*}
 H^0\left(\bigoplus_{i=1}^l \bar{K}^{(i)}_2(m_2)\right)
 \otimes{\mathcal A}(-m_2)
 \longrightarrow
 H^0\left(\bigoplus_{i=1}^l \bar{K}^{(i)}_1(m_1)\right)
 \otimes{\mathcal A}(-m_1) \\
 \longrightarrow
 H^0\left(\bigoplus_{i=1}^l E^{(i)}(m_0)\right)\otimes{\mathcal A}(-m_0)
 \longrightarrow \bigoplus_{i=1}^l E^{(i)} \longrightarrow 0.
\end{gather*}
Then we can see that the quiver representations determined by
$y$ and $y'$ are $S$-equivalent.
So we have $\phi(x)=\phi(y)=\phi(y')$.
Note that $G_s y'$ is a closed orbit in $(Y\cap Z^{ss}(\chi))_s$
by [\cite{king}, Proposition 3.2].
Thus the closure of the $G_s$-orbit of $x$ must contain $y'$.
Then, by Proposition \ref{openness-stability}, $x$ becomes a geometric point of $Y^{ss}_s$.
\end{proof}

\noindent
{\it Proof of Theorem \ref{thm:moduli-exists}.}
If we put
\[
 \overline{M_{\mathcal D}^{P,\mathcal L}}:=\phi(Y^{ss}),
\]
then we can see by Lemma \ref{open-lemma} that
$\overline{M_{\mathcal D}^{P,\mathcal L}}$ is an open subset of
$(Y\cap Z^{ss}(\chi))//G$.
We can see by a similar argument to that of [\cite{maruyama}, Proposition 7.3]
that there is a canonical morphism
$\Phi:\overline{{\mathcal M}_{\mathcal D}^{P,\mathcal L}}\rightarrow
\overline{M_{\mathcal D}^{P,\mathcal L}}$.
For two geometric points $x_1,x_2\in Y^{ss}$ over a geometric point $s$
of $S$,
$\phi(x_1)=\phi(x_2)$ if and only if the corresponding
representations of quiver are $S$-equivalent (\cite{king}),
that is, the corresponding objects of ${\mathcal D}_s$
are $S$-equivalent.
Thus for any algebraically closed field $k$ over $S$,
$\Phi(k):\overline{{\mathcal M}_{\mathcal D}^{P,\mathcal L}}(k)
\rightarrow \overline{M_{\mathcal D}^{P,\mathcal L}}(k)$
is bijective.
We can see by a standard argument that
$\overline{M_{\mathcal D}^{P,\mathcal L}}$ has the universal
property of the coarse moduli scheme.
If we put $M_{\mathcal D}^{P,\mathcal L}:=Y^s/G$,
then $M_{\mathcal D}^{P,\mathcal L}$ becomes an open subset of
$\overline{M_{\mathcal D}^{P,\mathcal L}}$
and we can easily see that $M_{\mathcal D}^{P,\mathcal L}$
is a coarse moduli scheme of
${\mathcal M}_{\mathcal D}^{P,\mathcal L}$.
So we have proved Theorem \ref{thm:moduli-exists}.
\hfill
$\square$

\begin{theorem}\label{projective}
 Assume that $S$ is of finite type over a universally Japanese ring $\Xi$.
 Then the moduli scheme $\overline{M_{\mathcal D}^{P,\mathcal L}}$
 is projective over $S$.
\end{theorem}

For the proof of Theorem \ref{projective},
the following lemma is essential.

\begin{lemma}\label{langton-thm}
 Let $R$ be a discrete valuation ring over $S$ with quotient field $K$
 and residue field $k$.
 Assume that $E$ is an object of ${\mathcal D}_K$
 which is ${\mathcal L}$-semistable.
 Then there is an object
 $\tilde{E}\in{\mathcal D}_R$ such that
 $\tilde{E}_K \cong E$
 and $\tilde{E}_k$ is ${\mathcal L}$-semistable.
\end{lemma}

\begin{proof}
The above $E$ corresponds to a coherent ${\mathcal A}_K$-module ${\mathcal E}$
and it suffices to show that there exists an $R$-flat coherent
${\mathcal A}_R$-module $\tilde{\mathcal E}$ such that
$\tilde{\mathcal E}\otimes_RK \cong {\mathcal E}$
and $\tilde{\mathcal E}\otimes k$ satisfies the semistability
condition given by the inequality in Remark \ref{stability-condition}.
For a sufficiently large integer $N$, we have
$H^i({\mathcal E}(N))=0$ for $i>0$ and
${\mathcal E}(N)$ is generated by global sections.
Then there is a surjection ${\mathcal A}_K(-N)^{\oplus r}
\rightarrow{\mathcal E}$
which determies a $K$-valued point $\eta$ of the Quot-scheme
$\Quot^P_{{\mathcal A}(-N)^{\oplus r}}$ for some numerical polynomial $P$,
where $r=\dim H^0({\mathcal E}(N))$.
Let ${\mathcal F}\subset{\mathcal A}(-N)^{\oplus r}$ be the universal
subsheaf and $Y$ be the maximal closed subscheme of
$\Quot^P_{{\mathcal A}(-N)^{\oplus r}}$ such that
${\mathcal A}\otimes {\mathcal F}_Y\rightarrow
{\mathcal A}(-N)^{\oplus r}_Y$
factors through ${\mathcal F}_Y$.  
Then $\eta$ is a $K$-valued point of $Y$ and extends to an $R$-valued point
$\xi$ of $Y$ because $Y$ is proper over $S$.
$\xi$ corresponds to an $R$-flat quotient coherent ${\mathcal A}_R$-module
${\mathcal E}'$ of ${\mathcal A}(-N)_R^{\oplus r}$
and we have ${\mathcal E}'\otimes_RK \cong {\mathcal E}$.
From the proof similar to that of Langton's theorem
(\cite{H-L}, Theorem 2.B.1),
we can obtain an $R$-flat coherent ${\mathcal A}_R$-module
$\tilde{\mathcal E}$ by taking succesive elementary transforms
of ${\mathcal E}'$ along $\mathbf{P}(V_1)\times\Spec k$
such that
$\tilde{\mathcal E}\otimes_RK\cong{\mathcal E}'\otimes_RK\cong{\mathcal E}$
and $\tilde{\mathcal E}\otimes k$ is semistable as
${\mathcal A}\otimes k$-module.
\end{proof}

Now we prove Theorem \ref{projective}.
By construction, the moduli scheme $\overline{M_{\mathcal D}^{P,\mathcal L}}$
is quasi-projective over $S$.
So it is sufficient to show that $\overline{M_{\mathcal D}^{P,\mathcal L}}$ is proper over $S$.
Let $R$ be a discrete valuation ring over $S$ with quotient field $K$
and let $\varphi: \Spec K \rightarrow \overline{M^{P,{\mathcal L}}_{\mathcal D}}$
be a morphism over $S$.
Then there is a finite extension field $K'$ of $K$
such that the composite
$\psi:\Spec K' \rightarrow \Spec K \stackrel{\varphi}\longrightarrow
\overline{M^{P,{\mathcal L}}_{\mathcal D}}$
is given by an ${\mathcal L}$-semistable object $E'$.
We can take a discrete valuation ring $R'$ with quotient field $K'$
such that $K\cap R'=R$.
Let $k'$ be the residue field of $R'$.
By Lemma \ref{langton-thm}, there exists an object $E$ of ${\mathcal D}_{R'}$
such that $E_{K'}\cong E'$ and $E_{k'}$ is ${\mathcal L}$-semistable.
Then $E$ gives a morphism
$\overline{\psi}:\Spec R'\rightarrow\overline{M^{P,{\mathcal L}}_{\mathcal D}}$
which is an extension of $\psi$.
We can easily see that $\overline{\psi}$ factors through $\Spec R$.
Thus $\overline{M^{P,{\mathcal L}}_{\mathcal D}}$ is proper over $S$
by the valuative criterion of properness.
\hfill
$\square$

\section{Examples}

In this section, we give several examples of moduli spaces of
stable objects determined by a strict ample sequence.

\begin{example}\label{example:stable-sheaf}\rm
Let $f:X\rightarrow S$ be a flat projective morphism of noetherian schemes
and let ${\mathcal O}_X(1)$ be an $S$-very ample line bundle on $X$ such that
$H^i({\mathcal O}_{X_s}(m))=0$ for $i>0$, $s\in S$ and $m>0$.
Consider the fibered triangulated category
${\mathcal D}_{X/S}$ defined by
$({\mathcal D}_{X/S})_U=D^b(\mathrm{Coh}(X_U/U))$ for $U\in(\sch/S)$.
Then ${\mathcal L}=\{ {\mathcal O}_X(-n) \}_{n\geq 0}$
becomes a strict ample sequence in ${\mathcal D}_{X/S}$.

\begin{proof}
Definition \ref{def-polarization} (1),(2),(3) are easy to verify.
Let us prove Definition \ref{def-polarization} (4).
Take any $U\in(\sch/S)$ and any object $E^{\bullet}\in({\mathcal D}_{X/S})_U$.
We may assume that $E^{\bullet}$ is given by a complex
\[
 \cdots \longrightarrow 0\longrightarrow 0 \longrightarrow
 E^{l_1}\stackrel{d^{l_1}}\longrightarrow E^{l_1+1} \stackrel{d^{l_1}+1}\longrightarrow
 \cdots \stackrel{d^{l_2-1}}\longrightarrow E^{l_2} \longrightarrow 0 \longrightarrow 0 \longrightarrow\cdots,
\]
where each $E^i$ is a coherent sheaf on $X_U$ flat over $U$.
By flattening stratification theorem, there is a stratification
$U=\coprod_{j=1}^m Y_j$ of $U$ by subschemes $Y_j$ such that
each $\coker(d^i)_{Y_j}=\coker(d^i_{Y_j})$ is flat over $Y_j$ for any $i$ and $j$.
Then we can see that $\im(d^i_{Y_j})$ and $\ker(d^i_{Y_j})$
are flat over $Y_j$ for any $i$ and $j$.
For any point $s\in U$, the sequence
\[
 0 \longrightarrow \im(d^{i-1}_{Y_j})\otimes k(s) \longrightarrow E^i\otimes k(s)
 \longrightarrow \coker(d^i_{Y_j})\otimes k(s) \longrightarrow 0 
\]
is exact because $\coker(d_{Y_j})$ is flat over $Y_j$.
Then the homomorphism
$\im(d^{i-1}_{Y_j})\otimes k(s) \longrightarrow \ker(d^i_{Y_j})\otimes k(s)$
is injective for any $s\in Y_j$.
Thus the cohomology sheaf
${\mathcal H}^i(E^{\bullet}_{Y_j}):=\ker(d^i_{Y_j})/\im(d^{i-1}_{Y_j})$
is flat over $Y_j$ for any $i$ and $j$.
We can take a positive integer $n_0$ such that for any $n\geq n_0$,
$R^p(f_{Y_j})_*(E^i_{Y_j}(n))=0$, $R^p(f_{Y_j})_*(\im(d^i_{Y_j})(n))=0$
and $R^p(f_{Y_j})_*(\ker d^i_{Y_j}(n))=0$
for any $p>0$ and any $i,j$.
Then we have $R^p(f_{Y_j})_*({\mathcal H}^i(E^{\bullet}_{Y_j}(n)))=0$
for any $p>0$, any $i,j$ and $n\geq n_0$.
From the spectral sequence
$R^p(f_{Y_j})_*({\mathcal H}^q(E^{\bullet}_{Y_J}(n)))\Rightarrow R^{p+q}(f_{Y_j})_*(E^{\bullet}_{Y_j}(n))$,
we have an isomorphism
$R^i(f_{Y_j})_*(E^{\bullet}_{Y_j}(n))\cong (f_{Y_j})_*({\mathcal H}^i(E^{\bullet}_{Y_j})(n))$
for any $i,j$ and $n\geq n_0$.
So we can see that $\mathbf{R}(f_{Y_j})_*(E^{\bullet}_{Y_j}(n))$ is quasi-isomorphic to the complex
\[
 \cdots\longrightarrow 0\longrightarrow (f_{Y_j})_*(E^{l_1}_{Y_j}(n)) \longrightarrow
 (f_{Y_j})_*(E^{l_1+1}_{Y_j}(n))\longrightarrow \cdots \longrightarrow
 (f_{Y_j})_*(E^{l_2}_{Y_j}(n)) \longrightarrow 0 \longrightarrow \cdots
\]
for any $i,j$ and $n\geq n_0$.
Note that there are canonical isomorphisms
\[
 \mathbf{H}^i(E^{\bullet}_s(n))\cong R^i(f_{Y_j})_*(E^{\bullet}_{Y_j}(n))\otimes k(s) \cong 
 (f_{Y_j})_*({\mathcal H}^i(E^{\bullet}_{Y_j})(n))\otimes k(s)\cong H^0(X_s,{\mathcal H}^i(E^{\bullet}_s)(n)).
\]
for any $i,j$, any $s\in Y_j$ and $n\geq n_0$.
If we take $n_0$ sufficiently larger, we may assume that the homomorphism
\[
 (f_{Y_j})^*(f_{Y_j})_*({\mathcal H}^i(E^{\bullet}_{Y_j}(n)))\longrightarrow {\mathcal H}^i(E^{\bullet}_{Y_j})(n)
\]
is surjective for any $n\geq n_0$ and any $i,j$.
Thus there exists a positive integer $N_0\gg n$ such that
\[
 (f_{Y_j})_*({\mathcal O}_{X_{Y_j}}(N-n))\otimes (f_{Y_j})_*({\mathcal H}^i(E^{\bullet}_{Y_j}(n)))
 \longrightarrow (f_{Y_j})_*({\mathcal H}^i(E^{\bullet}_{Y_j})(N))
\]
is surjective for any $N\geq N_0$ and any $i,j$.
So we obtain a commutative diagram
\[
 \begin{CD}
 (f_{Y_j})_*({\mathcal O}_{X_{Y_j}}(N-n))\otimes k(s)\otimes (f_{Y_j})_*({\mathcal H}^i({\mathcal E}^{\bullet}_{Y_j}(n)))\otimes k(s)
 @>>> (f_{Y_j})_*({\mathcal H}^i(E^{\bullet}_{Y_j}(N))\otimes k(s)   \\
 @V\cong VV @V\cong VV   \\
 H^0({\mathcal O}_{X_s}(N-n))\otimes \mathbf{H}^i(E^{\bullet}_s(n)) @>>> \mathbf{H}^i(E^{\bullet}_s(N)) \\
 @V \cong VV    @V \cong VV \\
 \Hom({\mathcal O}_{X_s}(-N),{\mathcal O}_{X_s}(-n))\otimes\Ext^i({\mathcal O}_{X_s}(-n),E^{\bullet}_s)
 @>>> \Ext^i({\mathcal O}_{X_s}(-N),E^{\bullet}_s).
 \end{CD}
\]
for any $i,j$, any $s\in Y_j$ and $N\geq N_0$.
Hence 
\[
 \Hom({\mathcal O}_{X_s}(-N),{\mathcal O}_{X_s}(-n))\otimes\Ext^i({\mathcal O}_{X_s}(-n),E^{\bullet}_s)
 \longrightarrow  \Ext^i({\mathcal O}_{X_s}(-N),E^{\bullet}_s)
\]
is surjective for any $s\in U$, any $i$ and $N\geq N_0$
and we have proved Definition \ref{def-polarization} (4).

Now we prove Definition \ref{def-polarization} (5).
Assume that an object $E\in({\mathcal D}_{X/S})_U$ and integers $i$, $n_0$ are given such that
$\Ext^i({\mathcal O}_{X_s}(-n),E^{\bullet}_s)=0$
for any $s\in U$ and $n\geq n_0$.
Replacing $n_0$ by a sufficiently large integer, we have
\[
 \Ext^i({\mathcal O}_{X_s}(-n),E^{\bullet}_s)\cong \mathbf{H}^i(E^{\bullet}_s(n))\cong H^0(X_s,{\mathcal H}^i(E^{\bullet}_s)(n))=0
\]
for any $s\in U$ and any $n\geq n_0$.
Then we have ${\mathcal H}^i(E^{\bullet}_s)=0$.
If $E^{\bullet}$ is given by 
\[
 E^{l_1}\stackrel{d^{l_1}}\longrightarrow E^{l_1+1}\stackrel{d^{l_1+1}}\longrightarrow
 \cdots\stackrel{d^{l_2-1}}\longrightarrow E^{l_2},
\]
such that each $E^j$ is flat over $U$, then the induced homomorphism
$\coker(d^{i-1})\otimes k(s) \rightarrow E^{i+1}\otimes k(s)$ is injective for any $s\in U$.
Then $\coker(d^i)$ is flat over $U$ and
$\coker(d^{i-1})\rightarrow E^{i+1}$ is injective.
Let $F^{\bullet}$ be the complex given by
\[
 \cdots\longrightarrow 0\longrightarrow \coker(d^i) \longrightarrow E^{i+2}
 \stackrel{d^{i+2}}\longrightarrow\cdots
 \stackrel{d^{l_2-1}}\longrightarrow E^{l_2} \longrightarrow 0 \longrightarrow \cdots.
\]
Then there is a canonical morphism
$u:E^{\bullet}\rightarrow F^{\bullet}$.
Note that
\[
 R^j\Hom_f({\mathcal O}_{X_U}(-n),E^{\bullet})=R^j(f_U)_*(E^{\bullet}(n))\cong
 (f_U)_*({\mathcal H}^j(E^{\bullet})(n))
\]
for $n\gg 0$.
So $u$ induces isomorphisms
\[
  R^j\Hom_f({\mathcal O}_{X_U}(-n),E^{\bullet})\stackrel{\sim}\longrightarrow
  (f_U)_*({\mathcal H}^j(E^{\bullet})(n))\stackrel{\sim}\longrightarrow
  (f_U)_*({\mathcal H}^j(F^{\bullet})(n))\stackrel{\sim}\longrightarrow
  R^j\Hom_f({\mathcal O}_{X_U}(-n),F^{\bullet})
  \]
 for $j>i$ and $n\gg 0$.
By definition we have $ R^j\Hom_f({\mathcal O}_{X_U}(-n),F^{\bullet})= (f_U)_*({\mathcal H}^j(F^{\bullet}(n)))=0$
for $j\leq i$ and $n\gg 0$.
Thus we have proved Definition \ref{def-polarization} (5).

Finally, let us prove Definition \ref{def-polarization} (6).
Let $E^{\bullet}$ and $F^{\bullet}$ be objects of $({\mathcal D}_{X/S})_U$.
Assume that $R^j(f_U)_*(E^{\bullet}(n))=0$ for $j\geq 0$ and $n\gg 0$ and that
$R^j(f_U)_*(F^{\bullet}(n))=0$ for $j<0$ and $n\gg 0$.
Since $R^j(f_U)_*(E^{\bullet}(n))\cong (f_U)_*({\mathcal H}^j(E^{\bullet})(n))$ for $n\gg 0$,
we have ${\mathcal H}^j(E^{\bullet})=0$ for $j\geq 0$.
Then $E^{\bullet}$ is quasi-isomorphic to the complex given by
\[
 \cdots\longrightarrow 0\longrightarrow E^{l_1}\stackrel{d^{l_1}_E}\longrightarrow E^{l_1+1}\longrightarrow
 \cdots\longrightarrow E^{-2} \longrightarrow \ker(d^{-1}_E) \longrightarrow 0\longrightarrow \cdots.
\]
On the other hand, we have ${\mathcal H}^j(F^{\bullet})=0$ for $j<0$, because
$R^j(f_U)_*(F^{\bullet}(n))\cong (f_U)_*({\mathcal H}^j(F^{\bullet})(n))$ for $n\gg 0$.
Then $F^{\bullet}$ is quasi-isomorphic to the complex given by
\[
 \cdots\longrightarrow 0 \longrightarrow \coker d^{-1}_F\longrightarrow F^1\stackrel{d^1_F}\longrightarrow
 \cdots\longrightarrow F^{m_2} \longrightarrow 0 \longrightarrow\cdots.
\]
We can take a complex
\[
 \cdots\longrightarrow 0\longrightarrow I^0\longrightarrow I^1\longrightarrow I^2\longrightarrow \cdots
\]
such that each $I^j$ is an injective sheaf on $X_U$ and that $I^{\bullet}$ is quasi-isomorphic to $F^{\bullet}$.
Then we have
$\Hom_{({\mathcal D}_{X/S})_U}(E^{\bullet},F^{\bullet})\cong H^0(\Hom^{\bullet}(E^{\bullet},I^{\bullet}))=0$.
So we have proved Definition \ref{def-polarization} (6).
\end{proof}

For an object $E\in({\mathcal D}_{X/S})_U$,
$\Ext^i({\mathcal O}_{X_s}(-n),E_s)=0$ for $n\gg 0$, $i\neq 0$ and $s\in U$
if and only if $E^{\bullet}$ is quasi-isomorphic to a coherent sheaf on $X_U$ flat over $U$.
Hence, for a numerical polynomial $P$,
the moduli space $M^{P,{\mathcal L}}_{{\mathcal D}_{X/S}}$
(resp.\ $\overline{M^{P,{\mathcal L}}_{{\mathcal D}_{X/S}}}$)
is just the usual moduli space of ${\mathcal O}_X(1)$-stable sheaves
(resp.\  moduli space of $S$-equivalence classes of
${\mathcal O}_X(1)$-semistable sheaves)
on $X$ over $S$.
\end{example}

\begin{example}\rm
Let $X$, $S$, ${\mathcal O}_X(1)$ and ${\mathcal D}_{X/S}$
be as in Example \ref{example:stable-sheaf}.
Take a vector bundle $G$ on $X$.
Replacing ${\mathcal O}_X(1)$ by some multiple,
${\mathcal L}_G=\{ G\otimes{\mathcal O}_X(-n) \}_{n\geq 0}$
also becomes a strict ample sequence in ${\mathcal D}_{X/S}$
and the moduli space
$M^{P,{\mathcal L}_G}_{{\mathcal D}_{X/S}}$
(resp.\ $\overline{M^{P,{\mathcal L}_G}_{{\mathcal D}_{X/S}}}$)
is the moduli space of $G$-twisted ${\mathcal O}_X(1)$-stable sheaves
(resp.\ moduli space of $S$-equivalence classes of $G$-twisted
${\mathcal O}_X(1)$-semistable sheaves) on $X$ over $S$. 
\end{example}

\begin{example}\label{ex:F-M}\rm
Let $X$, $Y$ be  projective schemes over an algebraically closed field $k$
and let ${\mathcal O}_X(1)$ be a very ample line bundle on $X$ such that
$H^i(X,{\mathcal O}_X(m))=0$ for $i>0$ and $m>0$.
Assume that a Fourier-Mukai transform
\begin{gather*}
 \Phi:D^b_c(X)\stackrel{\sim}\longrightarrow D^b_c(Y) \\
 E \mapsto \mathbf{R}(p_Y)_*(p_X^*(E)\otimes {\mathcal P})
\end{gather*}
with the kernel ${\mathcal P}\in D^b_c(X\times Y)$ is given.
Then $\Phi$ extends to an equivalence of fibered triangulated categories
\[
 \Phi:{\mathcal D}_{X/k} \stackrel{\sim}\longrightarrow {\mathcal D}_{Y/k}.
\]
Since ${\mathcal L}=\{ {\mathcal O}_X(-n) \}_{n\geq 0}$ is a strict ample
sequence in ${\mathcal D}_{X/k}$,
${\mathcal L}^{\Phi}=\{ \Phi({\mathcal O}_X(-n)) \}_{n\geq 0}$
is a strict ample sequence in ${\mathcal D}_{Y/k}$.
Moreover $\Phi$ determines an isomorphism
\[
 \Phi:M^{P,{\mathcal L}}_{{\mathcal D}_{X/k}} \stackrel{\sim}\longrightarrow
 M^{P,{\mathcal L}^{\Phi}}_{{\mathcal D}_{Y/k}}
\]
of the moduli space of stable sheaves on $X$ to the moduli space of
stable objects in $D^b_c(Y)$.
 
\end{example}

\begin{example}\rm
Let $G$ be a finite group and $X$ be a projective variety over $\mathbf{C}$
on which $G$ acts.
Take a $G$-linearized very ample line bundle ${\mathcal O}_X(1)$ on $X$
such that $H^i(X,{\mathcal O}_X(m))=0$ for $i>0$ and $m>0$.
Let $\rho_0,\rho_1,\ldots,\rho_s$ be the irreducible representations of $G$.
Consider the fibered triangulated category
${\mathcal D}^G_{X/\mathbf{C}}$ defined by
$({\mathcal D}^G_{X/\mathbf{C}})_U=D^G(\mathrm{Coh}(X_U/U))$,
for $U\in (\sch/\mathbf{C})$,
where $D^G(\mathrm{Coh}(X_U/U))$ is the full subcategory of
the derived category of bounded complexes of $G$-equivariant
coherent sheaves on $X_U$ consisting of the objects of finite
Tor-dimension over $U$.
For positive integers $r_0,r_1,\ldots,r_s$,
${\mathcal L}^G_{(r_0,\ldots,r_s)}=
\left\{ {\mathcal O}_X(-n)\otimes
(\rho_0^{\oplus r_0}\oplus\cdots\oplus\rho_s^{\oplus r_s})
\right\}_{n\geq 0}$
becomes a strict ample sequence in ${\mathcal D}^G_{X/\mathbf{C}}$.
The moduli space
$M^{P,{\mathcal L}^G_{(r_0,\ldots,r_s)}}_{{\mathcal D}^G_{X/\mathbf{C}}}$
is just the moduli space of $G$-equivariant sheaves ${\mathcal E}$ on $X$
satisfying the stability condition:
${\mathcal E}$ is of pure dimension $d=\deg P$ and
for any $G$-equivariant subsheaf $0\neq{\mathcal F}\subsetneq{\mathcal E}$,
the inequality
\[
 \genfrac{}{}{}{}{\Hom_G \left(
 \rho_0^{\oplus r_0}\oplus\cdots\oplus\rho_s^{\oplus r_s},
 H^0(X,{\mathcal F}\otimes{\mathcal O}_X(n)) \right)}
 {a_0({\mathcal F})}
 <
 \genfrac{}{}{}{}{\Hom_G \left(
 \rho_0^{\oplus r_0}\oplus\cdots\oplus\rho_s^{\oplus r_s},
 H^0(X,{\mathcal E}\otimes{\mathcal O}_X(n)) \right)}
 {a_0({\mathcal E})}
\]
holds for $n\gg 0$, where we define
\[
 \chi({\mathcal E}(m))=
 \sum_{i=0}^d a_i({\mathcal E})\begin{pmatrix} m+d-i \\ d-i \end{pmatrix}
 \quad \text{and}\quad
 \chi({\mathcal F}(m))=
 \sum_{i=0}^d a_i({\mathcal F})\begin{pmatrix} m+d-i \\ d-i \end{pmatrix}
\]
and so on.
\end{example}

\begin{example}\rm
Let $X$ be a projective variety over $\mathbf{C}$
and let ${\mathcal O}_X(1)$ be a very ample line bundle on $X$ such that
$H^i(X,{\mathcal O}_X(m))=0$ for $i>0$ and $m>0$.
For a torsion class $\alpha\in H^2(X,{\mathcal O}_X^{\times})$,
consider the fibered triangulated category ${\mathcal D}^{\alpha}_{X/\mathbf{C}}$
over $(\sch/\mathbf{C})$ defined by
$({\mathcal D}^{\alpha}_{X/{\mathbf{C}}})_U:=D^b(\mathrm{Coh}(X_U/U),\alpha_U)$,
where $D^b(\mathrm{Coh}(X_U/U),\alpha_U)$
is the derived category of bounded complexes of coherent
$\alpha_U$-twisted sheaves on $X\times U$ of finite Tor-dimension over $U$
and $\alpha_U$ is the image of $\alpha$ in
$H^2(X_U,{\mathcal O}_{X_U}^{\times})$.
For a locally free $\alpha$-twisted sheaf $G$ of finite rank on $X$,
${\mathcal L}^{\alpha}_G=\{ G\otimes{\mathcal O}_X(-n) \}_{n\geq 0}$
becomes a strict ample sequence in ${\mathcal D}^{\alpha}_{X/\mathbf{C}}$,
after replacing ${\mathcal O}_X(1)$ by some multiple.
The moduli space $M^{P,{\mathcal L}^{\alpha}_G}_{{\mathcal D}^{\alpha}_{X/\mathbf{C}}}$
is just the moduli space of $G$-twisted stable $\alpha$-twisted sheaves
on $X$ in the sense of \cite{yoshioka1}.
\end{example}

\noindent
{\bf Acknowledgments}
The author would like to thank Professor K\={o}ta Yoshioka
for valuable discussions and giving him a useful idea.
He also would like to thank Doctor Yukinobu Toda
for valuable discussions and comments.
The author is partly supported by Grant-in-Aid for
Scientific Research (Wakate-B 18740011)
Japan Society for the Promotion of Science.




%

\end{document}